\documentclass[12pt]{article}%
\usepackage{amsmath}
\usepackage{amsfonts}
\usepackage{amssymb}
\usepackage{hyperref}
\usepackage{tikz}
\usetikzlibrary{scopes, backgrounds}
\usepackage{graphicx}%
\providecommand{\U}[1]{\protect\rule{.1in}{.1in}}
\newtheorem{theorem}{Theorem}[section]

\newtheorem{claim}{Claim}

\newtheorem{corollary}[theorem]{Corollary}

\newtheorem{definition}{Definition}

\newtheorem{lemma}[theorem]{Lemma}

\newtheorem{problem}{Problem}
\newtheorem{proposition}[theorem]{Proposition}

\definecolor{plum}{RGB}{148,0,211}

\begin{document}

\title{\textbf{Lower~Boundary~Independent Broadcasts~in~Trees}}
\author{E. M. Marchessault\thanks{Partially funded by a Jamie Cassels Undergraduate
Research Award from the University of Victoria, Canada}\ \ and C. M.
Mynhardt\thanks{Funded by a Discovery Grant from the Natural Sciences and
Engineering Research Council of Canada}\\Department of Mathematics and Statistics\\University of Victoria\\PO BOX 1700 STN CSC\\Victoria, B.C.\\\textsc{Canada} V8W 2Y2\\{\small elisemarchessault@uvic.ca, kieka@uvic.ca}}
\maketitle

\begin{abstract}
A broadcast on a connected graph $G=(V,E)$ is a function $f:V\rightarrow
\{0,1,\dots,\operatorname{diam}(G)\}$ such that $f(v)\leq e(v)$ (the
eccentricity of $v$) for all $v\in V$ if $|V|\geq2$, and $f(v)=1$ if
$V=\{v\}$. The cost of $f$ is $\sigma(f)=\sum_{v\in V}f(v)$. Let $V_{f}%
^{+}=\{v\in V:f(v)>0\}$. A vertex $u$ hears $f$ from $v\in V_{f}^{+}$ if the
distance $d(u,v)\leq f(v)$. When $f$ is a broadcast such that every vertex $x$
that hears $f$ from more than one vertex in $V_{f}^{+}$ also satisfies
$d(x,u)\geq f(u)$ for all $u\in V_{f}^{+}$, we say that the broadcast only
overlaps in boundaries. A broadcast $f$ is boundary independent if it
overlaps only in boundaries. Denote by $i_{\operatorname{bn}}(G)$ the minimum
cost of a maximal boundary independent broadcast.

We obtain a characterization of maximal boundary independent broadcasts, show
that $i_{\operatorname{bn}}(T^{\prime})\leq i_{\operatorname{bn}}(T)$ for any
subtree $T^{\prime}$ of a tree $T$, and determine an upper bound for
$i_{\operatorname{bn}}(T)$ in terms of the broadcast domination number of $T$.
We show that this bound is sharp for an infinite class of trees.

\end{abstract}

\noindent\textbf{Keywords:\hspace{0.1in}}broadcast domination; broadcast
independence; boundary independence

\noindent\textbf{AMS Subject Classification Number 2010:\hspace{0.1in}}05C69

\section{Introduction}

If $X$ is an independent set of vertices of a graph $G$, then no edge of $G$
is incident with, or covered by, more than one vertex in $X$. Neilson \cite{N
thesis} used this property to generalize independent sets to independent
broadcasts, calling them boundary independent broadcasts (see Definition
\ref{bn-i}). Here we continue the investigation of the parameter
$i_{\operatorname{bn}}(G)$, the minimum cost of a maximal boundary independent
broadcast on $G$. We obtain conditions for a boundary independent broadcast to
be maximal boundary independent, prove that if $T^{\prime}$ is a subtree of a
tree $T$, then $i_{\operatorname{bn}}(T^{\prime})\leq i_{\operatorname{bn}%
}(T)$, and improve Neilson's upper bound for $i_{\operatorname{bn}}(G)$
(Theorem \ref{ThmBound}) in the case of trees, showing that it is sharp for
infinitely many trees.

We begin by giving definitions, background and known results in Section
\ref{Sec_Defs}. Section \ref{Sec_max_bn} contains the characterization of
maximal boundary independent broadcasts, as well as a corollary for later
reference. We prove the upper bound (Theorem \ref{Thm_ibn_bound}) in Section
\ref{Sec_bound}. A class of trees that demonstrates the sharpness of the bound
is given in Theorem \ref{Thm_Uniquely_radial}, which is proved in Section
\ref{Sec_unique}. This is preceded by preparatory results in Section
\ref{Sec_prep} and an informal outline of the proof in Section
\ref{Sec_outline}. We close by mentioning open problems in Section
\ref{Sec_open}.

We refer the reader to \cite{CLZ} for general graph theory concepts and
notation and to \cite{HHS} for domination related concepts. For a recent
survey of broadcasts in graphs, see the chapter by Henning, MacGillivray and
Yang \cite{HMY}. We denote domination number of a graph $G$ by $\gamma(G)$,
the independence number by $\alpha(G)$, and the minimum cardinality of a
maximal independent set (the \emph{independent domination number} of $G$)
by\emph{ }$i(G)$. If $P$ is a path in $G$, we denote the length of $P$ by
$\ell(P)$; a \emph{diametrical path} of $G$ is a path $P$ such that $\ell(P)=\operatorname{diam}(G)$ (the diameter of $G$).   

\section{Broadcast definitions and known results}

\label{Sec_Defs}A \emph{broadcast} on a connected graph $G=(V,E)$ is a
function $f:V\rightarrow\{0,1,\dots,\operatorname{diam}(G)\}$ such that
$f(v)\leq e(v)$ (the eccentricity of $v$) for all $v\in V$ if $|V|\geq2$, and
$f(v)=1$ if $V=\{v\}$. When $G$ is disconnected, we define a broadcast on $G$
as the union of broadcasts on its components. Let $V_{f}^{+}=\{v\in
V:f(v)>0\}$. A vertex in $V_{f}^{+}$ is called a \emph{broadcasting vertex}. A
vertex $u$ \emph{hears} $f$ from $v\in V_{f}^{+}$, and $v$ $f$%
-\emph{dominates} $u$, if the distance $d(u,v)\leq f(v)$. If $d(u,v)<f(v)$, we
also say that say that $v$ \emph{overdominates }$u$. A broadcast $f$ is
\emph{dominating} if every vertex of $G$ hears $f$ from at least one vertex in
$V_{f}^{+}$. The \emph{cost} of $f$ is $\sigma(f)=\sum_{v\in V}f(v)$, and the
\emph{broadcast domination number} (or simply the \emph{broadcast number}) of
$G$ is $\gamma_{b}(G)=\min\left\{  \sigma(f):f\text{ is a dominating broadcast
of }G\right\}  $. A dominating broadcast $f$ such that $\sigma(f)=\gamma
_{b}(G)$ is called a $\gamma_{b}$\emph{-broadcast}. Following \cite{MR}, for a
broadcast $f$ on $G$ and $v\in V_{f}^{+}$, we define the

\begin{itemize}
\item $f$-\emph{neighbourhood}$\ $of $v$ by $N_{f}(v)=\{u\in V:d(u,v)\leq
f(v)\}$,

\item $f$-\emph{boundary} of $v$ by $B_{f}(v)=\{u\in V:d(u,v)=f(v)\}$,

\item $f$-\emph{private neighbourhood}$\ $of $v$ by $\operatorname{PN}%
_{f}(v)=\{u\in N_{f}(v):u\notin N_{f}(w)$ for all$\ w\in V_{f}^{+}-\{v\}\}$,

\item $f$-\emph{private boundary}$\ $of $v$ by $\operatorname{PB}%
_{f}(v)=\{u\in N_{f}(v):u$ is not dominated by$\ (f-\{(v,f(v)\})\cup
\{(v,f(v)-1)\}$.
\end{itemize}

Note that if $f(u)=1$ and $u$ does not hear $f$ from any vertex $v\in
V_{f}^{+}-\{u\}$, then $u\in\operatorname{PB}_{f}(u)$ (but $u\notin B_{f}%
(u)$), and if $f(u)\geq2$,
then $\operatorname{PB}_{f}(u)=B_{f}(u)\cap\operatorname{PN}_{f}(u)$. If $f$
is a broadcast such that every vertex $x$ that hears more than one
broadcasting vertex also satisfies $d(x,u)\geq f(u)$ for all $u\in V_{f}^{+}$,
then the \emph{broadcast only overlaps in boundaries}. When $xy\in E(G)$ and
$x,y\in N_{f}(u)$ for some $u\in V_{f}^{+}$ such that at least one of $x$ and
$y$ does not belong to $B_{f}(u)$, we say that the edge $xy$ is \emph{covered}
in $f$ by $u$. When $xy$ is not covered by any $u\in V_{f}^{+}$, we say that
$xy$ is \emph{uncovered by~}$f$ or $f$-\emph{uncovered}. We denote the set of
$f$-uncovered edges by $U_{f}^{E}$.

\begin{figure}[ht]
    \centering
    \begin{tikzpicture}[main/.style = {draw, circle}, minimum size =0.25 cm, inner sep=0pt, node distance =0.75cm]
    
    \node[main,fill=white] (1) [left=2cm]{};
    \node[main,fill=white] (2) [right of =1] {};
    \node[main,fill=black] (3) [right of =2,label={[xshift=0.0cm, yshift=-0.75 cm]2}] {};
    \node[main,fill=white] (4) [right of =3] {};
    \node[main,fill=white] (5) [right of =4] {};
    \node[main,fill=white] (6) [right of =5] {};
    \node[main,fill=white] (7) [right of =6] {};
    \node[main,fill=black] (8) [right of =7,label={[xshift=0.0cm, yshift=-0.75 cm]2}] {};
    \node[main,fill=white] (9) [right of =8] {};
    \node[main,fill=white] (10) [right of =9] {};
    
    \node[main,fill=white] (11) [above of =2] {};
    \node[main,fill=white] (12) [above of =7] {};
    \node[main,fill=white] (13) [above of =8] {};
    \node[main,fill=white] (14) [above of =13] {};
    \node[main,fill=white] (15) [above of =9] {};
    
    \draw(1)--(2);
    \draw(2)--(3);
    \draw(3)--(4);
    \draw(4)--(5);
    \draw(5)--(6);
    \draw(6)--(7);
    \draw(7)--(8);
    \draw(8)--(9);
    \draw(9)--(10);
    \draw(2)--(11);
    \draw(7)--(12);
    \draw(13)--(8);
    \draw(14)--(13);
    \draw(15)--(9);

    {[on background layer]
    \draw[fill=orange,opacity=0.2] (-2.25,0) -- (1,0) -- (-0.625,1.7) -- cycle; 
    \draw[fill=orange,opacity=0.2] (1.5,0) -- (4.75,0) -- (3.125,1.7) -- cycle; 
    \draw[fill=orange,opacity=0.2] (5.9,0) -- (7.65,0) -- (6.755,0.95) -- cycle; 
    \draw[fill=orange,opacity=0.2] (8.1,0) -- (12.9,0) -- (10.5,2.5) -- cycle; 
    }

    \node[main,fill=white] (a) [right of =10, right = 0.5cm]{};
    \node[main,fill=black] (b) [right of =a,label={[xshift=0.0cm, yshift=-0.75 cm]1}] {};
    \node[main,fill=white] (c) [right of =b] {};
    \node[main,fill=white] (d) [right of =c] {};
    \node[main,fill=white] (e) [right of =d] {};
    \node[main,fill=white] (f) [right of =e] {};
    \node[main,fill=black] (g) [right of =f,label={[xshift=0.0cm, yshift=-0.75 cm]3}] {};
    \node[main,fill=white] (h) [right of =g] {};
    \node[main,fill=white] (i) [right of =h] {};
    \node[main,fill=white] (j) [right of =i] {};
    
    \node[main,fill=white] (k) [above of =b] {};
    \node[main,fill=white] (l) [above of =g] {};
    \node[main,fill=white] (m) [above of =h] {};
    \node[main,fill=white] (n) [above of =m] {};
    \node[main,fill=white] (o) [above of =i] {};
    
    \draw(a)--(b);
    \draw(b)--(c);
    \draw(c)--(d);
    \draw(d)--(e);
    \draw(e)--(f);
    \draw(f)--(g);
    \draw(g)--(h);
    \draw(h)--(i);
    \draw(i)--(j);
    \draw(k)--(b);
    \draw(l)--(g);
    \draw(m)--(h);
    \draw(m)--(n);
    \draw(o)--(i);
    
    \node at (1.25,0.25) {$e$};
    \node at (7.875,0.3) {$e'$};

    \end{tikzpicture}
    \caption{A tree with two maximum split-sets $\{e\}$ and $\{e'\}$; a central vertex and the radius of each subtree are shown}
    \label{Fig_split}
\end{figure}
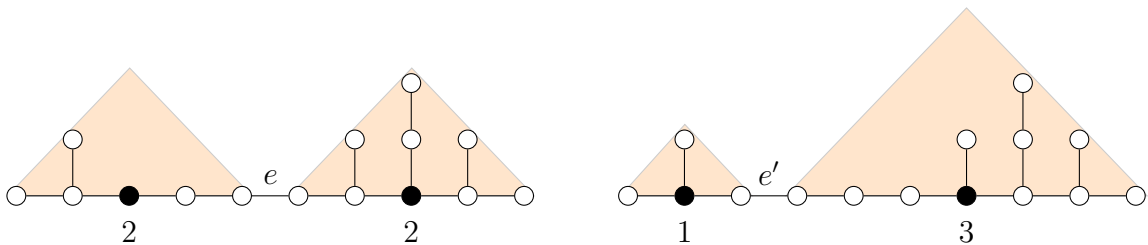

Erwin \cite{Erw thesis, Epaper} was the first to consider the broadcast
domination problem and noted that $\gamma_{b}(G)\leq\min\{\operatorname{rad}%
(G),\gamma(G)\}$. A natural question follows from this bound: which graphs $G$
satisfy $\gamma_{b}(G)=\operatorname{rad}(G)$? Such a graph $G$ is called
\emph{radial}; the problem of characterizing radial trees was addressed by
Dunbar, Erwin, Haynes, Hedetniemi and Hedetniemi \cite{DEHHH} and Dunbar,
Hedetniemi and Hedetniemi \cite{Dunbar}, and solved completely by Herke and
Mynhardt \cite{HM}. The characterization in \cite{HM} involves the concept of
a split-set, which is important for the present paper as well. Definition
\ref{Def_splitset} is illustrated in Figure \ref{Fig_split}.

\begin{definition}
\label{Def_splitset}\emph{Let }$T$\emph{ be a tree with diametrical path }%
$P$\emph{, and }$M\subseteq E(P)$\emph{ a set of cardinality }$m\geq1$\emph{.
Let }$T_{1},...,T_{m+1}$\emph{ be the components of }$T-M$\emph{, and }$P_{i}%
$\emph{ the subpath of }$P$\emph{ in~}$T_{i}$\emph{. \vspace{-0.1in}}

\begin{itemize}
\item \emph{We say that }$M$\emph{ is a }split-set\emph{ of }$T$\emph{ if each
}$P_{i}$\emph{ has even, positive length and is a diametrical path of~}$T_{i}%
$\emph{. \vspace{-0.1in}}

\item \emph{A split-set is a }maximum split-set\emph{ if it is a split-set of
}$T$\emph{ of maximum cardinality. \vspace{-0.1in}}

\item \emph{When }$M$\emph{ is a maximum split-set of }$T$\emph{, the
components of }$T-M$\emph{ are called the }radial subtrees \emph{of }%
$T$\emph{.}
\end{itemize}
\end{definition}

As shown in \cite{HM}, the cardinality of a maximum split-set is independent
of the chosen diametrical path of the tree. Herke and Mynhardt \cite{HM} used
Definition \ref{Def_splitset} to characterize radial trees and determine the
broadcast number of a given tree.

\begin{theorem}
\label{Thm_gammaT}\emph{\cite{HM}\hspace{0.1in}}$(i)\hspace{0.1in}$A tree $T$
is radial if and only if it does not possess a split-set.

$(ii)$ If $T$ is a radial, then $\gamma_{b}(T)=\left\lceil \frac{1}%
{2}\operatorname{diam}(T)\right\rceil =\operatorname{rad}(T)$. If $T$ is
nonradial and $M$ is a maximum split-set of $T$, then $\gamma_{b}(T)=\frac
{1}{2}(\operatorname{diam}(T)-|M|)$.
\end{theorem}

Dunbar, Hedetniemi and Hedetniemi \cite{Dunbar} showed that the broadcast
number of a tree is bounded below by that of any of its subtrees.

\begin{proposition}
\label{Prop_Subtree}\emph{\cite{Dunbar}}\hspace{0.1in}If $T^{\prime}$ is a
subtree of a tree $T$, then $\gamma_{b}(T^{\prime})\leq\gamma_{b}(T)$.
\end{proposition}

Erwin \cite{Erw thesis} defined a broadcast $f$ to be \emph{independent} if no
vertex $u\in V_{f}^{+}$ hears $f$ from any other vertex $v\in V_{f}^{+}$; that
is, broadcasting vertices only hear themselves. This version of broadcast
independence, which we call \emph{hearing independence}, was also considered
by, among others, Ahmane, Bouchemakh and Sopena \cite{ABS}, Bessy and
Rautenbach \cite{BR, BR2}, Bouchemakh and Zemir \cite{Bouch}, and Dunbar et
al.~\cite{DEHHH}. Consider the tree obtained from the star $K_{1,r},\ r\geq3$,
by subdividing each edge once. Broadcasting from each leaf with a strength of
$3$ gives a hearing independent broadcast of cost $3r$. Each non-pendant
vertex and edge hear the broadcast from $r$ vertices. Suppose, instead, that
we want as many edges as possible to hear (or be covered by) the broadcast,
but without signals overlapping on edges because such interference is
undesirable. To accomplish this, Mynhardt and Neilson \cite{MN} and Neilson
\cite{N thesis} pointed out that another way to generalize independent sets to
independent broadcasts is by ensuring that each edge belongs to the
$f$-neighbourhood of at most one broadcasting vertex; using this approach they
defined \emph{boundary independent broadcasts} as an alternative to hearing
independent broadcasts.

\begin{definition}
\label{bn-i}\emph{\cite{MN, N thesis}\hspace{0.1in}Let }$f$\emph{ be a
broadcast on a graph }$G$\emph{. We say that \vspace{-0.1in}}

\begin{itemize}
\item $f$\emph{ is }boundary independent, \emph{abbreviated}
bn-independent\emph{, if it overlaps only in boundaries,\vspace{-0.1in}}

\item \emph{and }maximal bn-independent\emph{ if, in addition, there is no
bn-independent broadcast }$g$\emph{ on }$G$\emph{ such that }$g>f$.\emph{
\vspace{-0.1in}}

\item \emph{We denote the minimum weight of a maximal bn-independent broadcast
}$f$\emph{ on }$G$\emph{ by }$i_{\operatorname{bn}}(G)$\emph{, and call }%
$f$\emph{ an }$i_{\operatorname{bn}}$-broadcast\emph{.}
\end{itemize}
\end{definition}

Erwin similarly denoted the minimum weight of a maximal hearing independent
broadcast $f$ on $G$ by $i_{b}(G)$\emph{.} Neilson \cite{N thesis} and
Mynhardt and Neilson \cite{MN} characterized bn-independent broadcasts that
are maximal bn-independent.

\begin{proposition}
\label{Thm_Minimal}\emph{\cite{MN, N thesis}}\hspace{0.1in}A bn-independent
broadcast $f$ is maximal bn-independent if and only if it is dominating, and
either $V_{f}^{+}=\{v\}$ or $B_{f}(v)-\operatorname{PB}_{f}(v)\neq\varnothing$
for each $v\in V_{f}^{+}$.
\end{proposition}

\begin{figure}[ht]
    \centering
    \begin{tikzpicture}[main/.style = {draw, circle}, minimum size =0.25 cm, inner sep=0pt, node distance =0.75cm]
    
    \node[main,fill=red] (1) [left=2cm]{};
    \node[main,fill=red] (2) [right of =1] {};
    \node[main,fill=red,ultra thick] (3) [right of =2,label={[xshift=0.0cm, yshift=-0.75 cm]\textcolor{red}{2}}] {};
    \node[main,fill=red] (4) [right of =3] {};
    \node[rectangle,draw,fill=plum] (5) [right of =4] {};
    \node[main,fill=blue] (6) [right of =5] {};
    \node[main,fill=blue] (7) [right of =6] {};
    \node[main,fill=blue,ultra thick] (8) [right of =7,label={[xshift=0.0cm, yshift=-0.75 cm]\textcolor{blue}{3}}] {};
    \node[main,fill=blue] (9) [right of =8] {};
    \node[main,fill=blue] (10) [right of =9] {};
    
    \node[main,fill=red] (11) [above of =2] {};
    \node[main,fill=blue] (12) [above of =7] {};
    \node[main,fill=blue] (13) [above of =8] {};
    \node[main,fill=blue] (14) [above of =13] {};
    \node[main,fill=blue] (15) [above of =9] {};
    
    \draw(1)--(2);
    \draw(2)--(3);
    \draw(3)--(4);
    \draw(4)--(5);
    \draw(5)--(6);
    \draw(6)--(7);
    \draw(7)--(8);
    \draw(8)--(9);
    \draw(9)--(10);
    \draw(2)--(11);
    \draw(7)--(12);
    \draw(13)--(8);
    \draw(14)--(13);
    \draw(15)--(9);

    \node[main,fill=red] (a) [right of =10, right = 0.5cm]{};
    \node[main,fill=red, ultra thick] (b) [right of =a,label={[xshift=0.0cm, yshift=-0.75 cm]\textcolor{red}{1}}] {};
    \node[rectangle,draw,fill=orange] (c) [right of =b] {};
    \node[main,fill=yellow,ultra thick] (d) [right of =c,label={[xshift=0.0cm, yshift=-0.75 cm]1}] {};
    \node[rectangle,draw,fill=green] (e) [right of =d] {};
    \node[main,fill=blue] (f) [right of =e] {};
    \node[main,fill=blue] (g) [right of =f] {};
    \node[main,fill=blue,ultra thick] (h) [right of =g,label={[xshift=0.0cm, yshift=-0.75 cm]\textcolor{blue}{3}}] {};
    \node[main,fill=blue] (i) [right of =h] {};
    \node[main,fill=blue] (j) [right of =i] {};
    
    \node[main,fill=red] (k) [above of =b] {};
    \node[main,fill=blue] (l) [above of =g] {};
    \node[main,fill=blue] (m) [above of =h] {};
    \node[main,fill=blue] (n) [above of =m] {};
    \node[main,fill=blue] (o) [above of =i] {};
    
    \draw(a)--(b);
    \draw(b)--(c);
    \draw(c)--(d);
    \draw(d)--(e);
    \draw(e)--(f);
    \draw(f)--(g);
    \draw(g)--(h);
    \draw(h)--(i);
    \draw(i)--(j);
    \draw(k)--(b);
    \draw(l)--(g);
    \draw(m)--(h);
    \draw(m)--(n);
    \draw(o)--(i);

    \end{tikzpicture}
    \caption{Two maximal bn-independent broadcasts $f$ and $g$ on the tree in Figure \ref{Fig_split}}
    \label{Fig_bn_ind}
\end{figure}
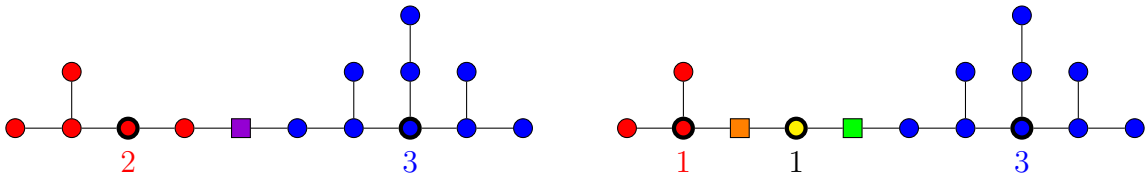


Two maximal bn-independent broadcasts $f$ and $g$ on the tree in Figure
\ref{Fig_split} are shown in Figure~\ref{Fig_bn_ind}. The broadcasting
vertices are shown in bold-outlined primary colours, the shared boundaries as
squares in the corresponding secondary colours, and the rest of the $f$- and
$g$-neighbourhoods in the same primary colours as the broadcasting vertices.

Proposition \ref{Thm_Minimal} implies that $\gamma_{b}(G)\leq
i_{\operatorname{bn}}(G)$ for every graph $G$. Moreover, broadcasting with a
strength of $\operatorname{rad}(G)$ from a central vertex of $G$ produces a
maximal bn-independent broadcast, hence $i_{\operatorname{bn}}(G)\leq
\operatorname{rad}(G)$ for all graphs $G$. In particular, if $G$ is radial,
then $\gamma_{b}(G)=i_{\operatorname{bn}}(G)$. On the other hand,
$i_{\operatorname{bn}}(G)$ is not comparable with the independent domination
number $i(G)$. For example, $i(P_{6})=2$ and $i_{\operatorname{bn}}(P_{6})=3$
\cite[Theorem 3.1.4]{N thesis}, whereas $i(K_{3,3})=3$ and
$i_{\operatorname{bn}}(K_{3,3})=\operatorname{rad}(K_{3,3})=2.$

Neilson \cite{N thesis} bounded $i_{\operatorname{bn}}(G)$ in terms of
$\gamma_{b}(G)$ as follows.

\begin{theorem}
\emph{\cite[Corollary 2.2.10]{N thesis}}\label{ThmBound}\hspace{0.1in}For any
graph $G$, $i_{\operatorname{bn}}(G)\leq\left\lceil \frac{4\gamma_{b}(G)}%
{3}\right\rceil $.
\end{theorem}

Our goal is to improve the bound in Theorem \ref{ThmBound} for trees as
follows, and to show that the bound is achieved by certain trees whose maximum
split-sets have cardinality $1$ or $2$.

\begin{theorem}
\label{Thm_ibn_bound}\emph{\hspace{0.1in}}Let $T$ be a non-radial tree having
a maximum split-set $M$. Then $i_{\operatorname{bn}}(T)\leq\gamma
_{b}(T)+\left\lceil \frac{|M|+1}{3}\right\rceil $.
\end{theorem}

\begin{theorem}
\label{Thm_Uniquely_radial}Let $T$ be a tree whose maximum split-set(s) have
cardinality $m\in\{1,2\}$. Suppose that for each maximum split-set $M$, every
component of $T-M$ is uniquely radial. Then $i_{\operatorname{bn}}%
(T)=\gamma_{b}(T)+1$.
\end{theorem}

\section{Maximal bn-independent broadcasts}

\label{Sec_max_bn}Before turning our attention to trees, we provide a
characterization, different from the one in Proposition \ref{Thm_Minimal}, of
maximal bn-independent broadcasts having at least two broadcasting vertices.
This characterization is of general interest and it and its corollary will be
used frequently in the proofs that follow. 

\begin{proposition}
\label{Prop_2bv}Let $f$ be a bn-independent broadcast on a connected graph $G$
such that $|V_{f}^{+}|\geq2$. Then $f$ is maximal bn-independent if and only
if each component of $G-U_{f}^{E}$ contains at least two broadcasting vertices.
\end{proposition}

\noindent\textbf{Proof}.\hspace{0.1in}Assume that $f$ is maximal
bn-independent. By Proposition~\ref{Thm_Minimal}, $f$ is dominating. Suppose
that some component $H$ of $G-U_{f}^{E}$ contains fewer than two broadcasting
vertices. Since $f$ is dominating, $|V_{f}^{+}\cap V(H)|=1$; say $V_{f}%
^{+}\cap V(H)=\{v\}$. Since $G$ is connected, $v$ is not isolated in $G$, and, since $v$ is broadcasting, all edges incident with $v$ also belong to $H$. Therefore $|V(H)|\geq2$. Since $H$ is a component of $G-U_{f}^{E}$, $B_{f}%
(v)\neq\varnothing$ and $B_{f}(v)\subseteq\operatorname{PB}_{f}(v)$. But then
$B_{f}(v)-\operatorname{PB}_{f}(v)=\varnothing$, contradicting
Proposition~\ref{Thm_Minimal}.

Conversely, assume that each component of $G-U_{f}^{E}$ contains at least two
broadcasting vertices. Since $G-U_{f}^{E}$ is a spanning subgraph of $G$, $f$
is dominating. Let $v$ be any vertex of $G$. First assume that $f(v)=0$. Since
$f$ is dominating, it is clear that increasing the strength of the broadcast
from $v$ results in a broadcast that is not bn-independent. Assume therefore
that $v\in V_{f}^{+}$ and say $v$ belongs to the component $H$ of $G-U_{f}%
^{E}$. Let $h$ be the restriction of $f$ to $H$; clearly, $h$ is
bn-independent. Since $h$ covers all edges of $H$ and $|V_{h}^{+}|\geq2$,
there exists a vertex $u\in V_{h}^{+}-\{v\}$ such that some vertex $w$ of $H$
hears the broadcast from $u$ as well as from $v$. Since $h$ is bn-independent,
$w\in B_{h}(u)\cap B_{h}(v)$. Thus $w\notin\operatorname{PB}_{h}(v)$, so
$B_{h}(v)-\operatorname{PB}_{h}(v)\neq\varnothing$. This implies that
$B_{f}(v)-\operatorname{PB}_{f}(v)\neq\varnothing$. The result follows from
Proposition~\ref{Thm_Minimal}.~$\blacksquare$

\begin{corollary}
\label{Cor_2bv}Let $f$ be a maximal bn-independent broadcast on a connected
graph $G$ with $|V_{f}^{+}|\geq2$ and let $X\subseteq U_{f}^{E}$. For each
component $H$ of $G-X$, denote the restriction of $f$ to $H$ by $f_{H}$. Then
$f_{H}$ is a maximal bn-independent broadcast on $H$. 
\end{corollary}

\noindent\textbf{Proof}.\hspace{0.1in}By Proposition \ref{Prop_2bv}, each
component of $G-U_{f}^{E}$ has at least two broadcasting vertices. Since
$X\subseteq U_{f}^{E}$, each component of $G-X$, including $H$, therefore has
at least two broadcasting vertices.

Let $e=xy$ be an $f_{H}$-uncovered edge of $H$. Then $e$ is also an
$f$-uncovered edge of $G$, otherwise there exists a vertex $v\in V_{f}%
^{+}-V(H)$ that broadcasts to $x$ and $y$. But then $G-U_{f}^{E}$ contains a
path from $v$ to $x$, hence so does $G-X$. Since $x\in V(H)$ and $v\notin
V(H)$, this contradicts $H$ being a component of $G-X$. Therefore each $f_{H}%
$-uncovered edge of $H$ is an $f$-uncovered edge of $G$. It follows that each
component of $H-U_{f_{H}}^{E}$ is a component of $G-U_{f}^{E}$, so each
component of $H-U_{f_{H}}^{E}$ has at least two broadcasting vertices. By
Proposition \ref{Prop_2bv}, $f_{H}$ is maximal bn-independent.~$\blacksquare$

\section{Bn-independent broadcasts on trees}

Our focus in this section is the proofs of Theorem \ref{Thm_ibn_bound} in
Section \ref{Sec_bound} and Theorem \ref{Thm_Uniquely_radial} in Section
\ref{Sec_unique}. The latter proof requires substantial preparation, which is
given in Section \ref{Sec_prep}. Some of these results, such as Lemma
\ref{Lem_edge_cover} and Theorem \ref{Prop_subtree_bn}, are also of wider interest.

\subsection{Proof of Theorem \ref{Thm_ibn_bound}}

\label{Sec_bound}We prove Theorem \ref{Thm_ibn_bound} using the same approach
as Neilson \cite[Proposition 2.2.9]{N thesis}, as well as the well-known fact
that $i(P_{n})=\left\lceil n/3\right\rceil $ \cite[Theorem 6.1]{HHS}. We
restate the theorem for convenience.

\bigskip

\noindent\textbf{Theorem \ref{Thm_ibn_bound}}\hspace{0.1in}\emph{Let }%
$T$\emph{ be a non-radial tree having a maximum split-set }$M$\emph{. Then
}$i_{\operatorname{bn}}(T)\leq\gamma_{b}(T)+\left\lceil \frac{|M|+1}%
{3}\right\rceil $\emph{.\bigskip}

\noindent\textbf{Proof}.\hspace{0.1in}Say $m=|M|$ and consider the $m+1$
radial subtrees $T_{1},...,T_{m+1}$ of $T-M$. By Definition \ref{Def_splitset}, each $T_{i}$ has even diameter, hence is a central tree. Let $v_{i}$ be the unique
central vertex of $T_{i}$. We know that $\gamma_{b}(T_{i})=\operatorname{rad}%
(T_{i})$ for each $i$. Theorem \ref{Thm_gammaT}$(ii)$ implies that $\gamma
_{b}(T)=\sum_{i=1}^{m+1}\operatorname{rad}(T_{i})$. Hence the broadcast $f$ on
$T$ defined by $f(v_{i})=\operatorname{rad}(T_{i})$ for $i=1,...,m+1$ and
$f(v)=0$ otherwise is a bn-independent $\gamma_{b}$-broadcast on $T$. Observe
that $M=U_{f}^{E}$ and that each component of $T-U_{f}^{E}$ contains exactly
one broadcasting vertex. By Proposition \ref{Prop_2bv}, $f$ is not maximal bn-independent.

Let $P\cong P_{m+1}$ be the path whose vertices and edges are the trees
$T_{i}$ and the edges in $M$, respectively. Since $i(P)=\left\lceil
(m+1)/3\right\rceil $, there exists an independent dominating set $X$ of $P$
with $|X|=\left\lceil (m+1)/3\right\rceil $. Let $\mathcal{X}$ denote the
collection of trees $T_{i}$ that correspond to a vertex in $X$, and
$\mathcal{Y}$ the rest of the trees $T_{i}$. Define the broadcast $g$ on $T$
by $g(v_{i})=f(v_{i})+1$ if $T_{i}\in\mathcal{X}$, and $g(v)=f(v)$ otherwise.
Then $\sigma(g)=\sigma(f)+\left\lceil (m+1)/3\right\rceil $.

The independence of $X$ ensures that each edge in $M$ (and hence each edge of
$T$) is covered by at most one vertex in $V_{g}^{+}$, while the fact that $X$
is dominating ensures that each $T_{j}\in\mathcal{Y}$ belongs to the same
component of $T-U_{g}$ as some $T_{i}\in\mathcal{X}$. Conversely, since $f$ is
a dominating broadcast, each $T_{i}\in\mathcal{X}$ belongs to the same
component of $T-U_{g}$ as some $T_{j}\in\mathcal{Y}$. Hence each component of
$T-U_{g}$ contains at least two vertices belonging to $V_{g}^{+}$. By
Proposition \ref{Prop_2bv}, $g$ is maximal bn-independent. Therefore
$i_{\operatorname{bn}}(T)\leq\sigma(g)=\gamma_{b}(T)+\left\lceil
(m+1)/3\right\rceil $\emph{.~}$\blacksquare$

\bigskip

\begin{figure}[ht]
    \centering
    
    \begin{tikzpicture}[main/.style = {draw, circle}, minimum size =0.25 cm, inner sep=0pt, node distance =0.75cm]
    
    \node[main,fill=red] (1) {};
    \node[very thick,main,fill=red] (2) [right of =1,label={[xshift=0.0cm, yshift=-0.75 cm]\textcolor{red}{1}}] {};
    \node[rectangle,draw,fill=plum] (3) [right of =2] {};
    \node[very thick,main,fill=blue] (4) [right of =3,label={[xshift=0.0cm, yshift=-0.75 cm]\textcolor{blue}{1}}] {};
    \node[main,fill=blue] (5) [right of =4] {};
    \node[main,fill=red] (6) [right of =5] {};
    \node[main,fill=red,very thick] (7) [right of =6,label={[xshift=0.0cm, yshift=-0.75 cm]\textcolor{red}{1}}] {};
    \node[fill=plum,rectangle,draw] (8) [right of =7] {};
    \node[main,fill=blue,very thick] (9) [right of =8,label={[xshift=0.0cm, yshift=-0.75 cm]\textcolor{blue}{1}}] {};
    \node[main,fill=blue] (10) [right of =9] {};
    \node[main,fill=red] (11) [right of =10] {};
    \node[main,fill=red,very thick] (12) [right of =11,label={[xshift=0.0cm, yshift=-0.75 cm]\textcolor{red}{1}}] {};
    \node[fill=plum,rectangle,draw] (13) [right of =12] {};
    \node[main,fill=blue,very thick] (14) [right of =13,label={[xshift=0.0cm, yshift=-0.75 cm]\textcolor{blue}{1}}] {};
    \node[main,fill=blue] (15) [right of =14] {};
    
    \node[main,fill=red] (16) [above of =2] {};
    \node[main,fill=blue] (17) [above of =4] {};
    \node[main,fill=red] (18) [above of =7] {};
    \node[main,fill=blue] (19) [above of =9] {};
    \node[main,fill=red] (20) [above of =12] {};
    \node[main,fill=blue] (21) [above of =14] {};
    
    \draw(1)--(2);
    \draw(2)--(3);
    \draw(3)--(4);
    \draw(4)--(5);
    \draw[ultra thick](5)--(6);
    \draw(6)--(7);
    \draw(7)--(8);
    \draw(8)--(9);
    \draw(9)--(10);
    \draw[ultra thick](10)--(11);
    \draw(11)--(12);
    \draw(12)--(13);
    \draw(13)--(14);
    \draw(14)--(15);
    \draw(16)--(2);
    \draw(4)--(17);
    \draw(18)--(7);
    \draw(19)--(9);
    \draw(20)--(12);
    \draw(21)--(14);
    
    {[on background layer]
    \draw[fill=orange,opacity=0.2] (-0.1,0) -- (3.1,0) -- (1.5,1.8) -- cycle; 
    \draw[fill=orange,opacity=0.2] (3.7,0) -- (6.8,0) -- (5.25,1.8) -- cycle; 
    \draw[fill=orange,opacity=0.2] (7.4,0) -- (10.6,0) -- (9,1.8) -- cycle;
    
    \draw[fill=orange,opacity=0.2] (-0.1,-3.9) -- (4.6,-3.9) -- (2.25,-1.4) -- cycle;
    
    \draw[fill=orange,opacity=0.2] (5.2,-3.9) -- (9.85,-3.9) -- (7.5,-1.4) -- cycle;
    }

    \node[] at (-0.75,0.75) {$T_1$};
    \node[] at (3.4,0.25) {$e$};
    \node[] at (7.1,0.3) {$e'$};
    \node[] at (4.9, -3.65) {$e$};

    \node[main,fill=red] (a) [below of =1,below=3cm] {};
    \node[main,fill=red] (b) [right of =a] {};
    \node[main,fill=red] (c) [right of =b] {};
    \node[main,very thick,fill=red] (d) [right of =c,label={[xshift=0.0cm, yshift=-0.75 cm]\textcolor{red}{3}}] {};
    \node[main,fill=red] (e) [right of =d] {};
    \node[main,fill=red] (f) [right of =e] {};
    \node[main,fill=red] (g) [right of =f] {};
    \node[main,fill=yellow] (h) [right of =g] {};
    \node[main,fill=yellow, very thick] (i) [right of =h,label={[xshift=0.0cm, yshift=-0.75 cm]1}] {};
    \node[main,fill=yellow] (j) [right of =i] {};
    \node[main,fill=red] (k) [right of =j]{};
    \node[main,fill=blue] (l) [right of =k] {};
    \node[main,fill=blue, very thick] (m) [right of =l,label={[xshift=0.0cm, yshift=-0.75 cm]\textcolor{blue}{1}}] {};
    \node[main,fill=blue] (n) [right of =m] {};
    
    \node[main,fill=red] (o) [above of =d] {};
    \node[main,fill=red] (p) [above of =o] {};
    \node[main,fill=red,very thick] (q) [above of =k,label={[xshift=0.35cm, yshift=-0.25 cm]\textcolor{red}{1}}] {};
    \node[main,fill=red] (r) [above of =q] {};
    
    \draw(a)--(b);
    \draw(b)--(c);
    \draw(c)--(d);
    \draw(d)--(e);
    \draw(e)--(f);
    \draw(f)--(g);
    \draw[ultra thick](g)--(h);
    \draw(h)--(i);
    \draw(i)--(j);
    \draw[ultra thick](j)--(k);
    \draw[ultra thick](k)--(l);
    \draw(l)--(m);
    \draw(m)--(n);
    
    \draw(o)--(d);
    \draw(p)--(o);
    \draw(q)--(k);
    \draw(r)--(q);
    
    \node[] at (-0.75,-3.25) {$T_2$};
    
    \end{tikzpicture}
    \caption{Trees $T_1$ and $T_2$ whose radial subtrees are not uniquely radial, with a $\gamma_b$- and $i_{\operatorname{bn}}$-broadcast on $T_1$, and a $\gamma_b$-broadcast on $T_2$ that is not maximal bn-independent}
    \label{Fig_Bound}
\end{figure}
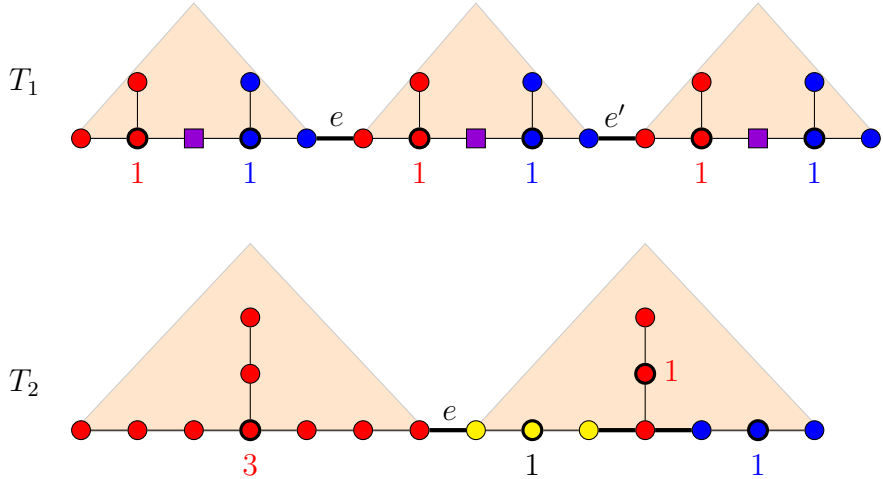

Observe that the bound in Theorem \ref{Thm_ibn_bound} is better than the one
in Theorem \ref{ThmBound} for all values of $|M|$ when $\gamma_{b}(T)\geq|M|+3$,
that is, when $\sum_{i=1}^{|M|+1}\operatorname{rad}(T_{i})\geq|M|+3$. For example, it is shown below that for the tree $T_{1}$ in Figure \ref{Fig_Bound}, $i_{\operatorname{bn}%
}(T_{1})=\gamma_{b}(T_{1})=6$. The bound in Theorem \ref{ThmBound} gives $i_{\operatorname{bn}%
}(T_{1})\leq8$ whereas the bound in Theorem \ref{Thm_ibn_bound} gives $i_{\operatorname{bn}%
}(T_{1})\leq7$. 

As examples of trees whose radial subtrees are not uniquely radial, consider
the trees in Figure \ref{Fig_Bound}. The tree $T_{1}$ has a unique maximum
split-set $\{e,e^{\prime}\}$. Each of its radial subtrees has a dominating
broadcast of cost $2$ from a central vertex, as well as the broadcast shown.
The uncovered edges are shown in cyan. By Proposition \ref{Prop_2bv}, the
broadcast on $T_{1}$ is maximal bn-independent. Hence $i_{\operatorname{bn}%
}(T_{1})=\gamma_{b}(T_{1})=6$. The tree $T_{2}$ has a unique maximum split-set
$\{e\}$. By Proposition \ref{Prop_2bv}, the broadcast on $T_{2}$ is not
maximal bn-independent. It is easy to check that the two broadcasts shown are
the only possible $\gamma_{b}$-broadcasts on the (isomorphic) radial subtrees
of $T_{2}$, and no combination of them produces a maximal bn-independent
broadcast on $T_{2}$. By Theorem \ref{Thm_ibn_bound}, $i_{\operatorname{bn}%
}(T_{2})=\gamma_{b}(T_{2})+1=7$.

\subsection{Preparatory results for the proof of Theorem
\ref{Thm_Uniquely_radial}}

\label{Sec_prep}We begin our preparation for the proof of Theorem
\ref{Thm_Uniquely_radial} with a lemma that gives a lower bound on the cost of
a broadcast $f$ on a tree $T$ in terms of the number of edges of a path $P$
that are covered by $f$. Each vertex $v\in V_{f}^{+}$ covers at most $2f(v)$
edges on $P$. Therefore, if $f$ covers $q$ edges on $P$, then $\sigma
(f)\geq\left\lceil \frac{q}{2}\right\rceil $. In fact, a stronger result
holds, as we show in Lemma \ref{Lem_edge_cover} below.

For any vertex $v$ of $T$, let $d(v,P)$ denote the distance from $v$ to the
vertex on $P$ nearest to $v$; note that $d(v,P)=0$ if and only if $v\in V(P)$.
With respect to $P$, the set $V_{f}^{+}$ of broadcasting vertices of $f$ are
of three types:\ those that lie on $P$, those that do not lie on $P$ but
broadcast to at least one vertex of $P$, and those that do not broadcast to
any vertex of $P$ at all. Denote the first and second types by
$\operatorname{Touch}(P)$, that is, $\operatorname{Touch}(P)=\{v\in V_{f}%
^{+}:N_{f}(v)\cap V(P)\neq\varnothing\}$, and the third type by
$\operatorname{Off}(P)$, that is, $\operatorname{Off}(P)=\{v\in V_{f}%
^{+}:N_{f}(v)\cap V(P)=\varnothing\}$. The next lemma holds for any broadcast
on a tree, not only bn-independent broadcasts.

\begin{lemma}
\label{Lem_edge_cover}Let $P$ be a path of a tree $T$ and let $f$ be a
broadcast on $T$. Suppose
\[
\sigma_{1}=\sum_{v\in\operatorname{Touch}(P)}d(v,P)\ \text{and}\ \sigma
_{2}=\sum_{v\in\operatorname{Off}(P)}f(v).
\]
Then

\begin{enumerate}
\item[$(i)$] $f$ covers at most $2(\sum_{v\in\operatorname{Touch}%
(P)}f(v)-\sigma_{1})$ edges of $P$;

\item[$(ii)$] if $f$ covers $b$ edges of $P$, then $\sigma(f)\geq\left\lceil
\frac{b}{2}\right\rceil +\sigma_{1}+\sigma_{2}$.
\end{enumerate}
\end{lemma}

\noindent\textbf{Proof.\hspace{0.1in}}Consider any vertex $v\in
\operatorname{Touch}(P)$ and let $u$ be the unique vertex of $P$ at distance
$d(v,P)$ from $v$. Since $u$ hears the broadcast from $v$, $f(v)$ covers the
same number of edges of $P$ as a broadcast with strength $f(v)-d(v,P)$ from
$u$ does, namely at most $2(f(v)-d(v,P))$. Summing over $\operatorname{Touch}%
(P)$, we obtain $(i)$.\emph{ }If $f$ covers $b$ edges of $P$, then $(i)$
implies that $b\leq2(\sum_{v\in\operatorname{Touch}(P)}f(v)-\sigma_{1})$, from
which it follows that $\sum_{v\in\operatorname{Touch}(P)}f(v)\geq\left\lceil
\frac{b}{2}\right\rceil +\sigma_{1}$. Since the vertices in
$\operatorname{Off}(P)$ do not broadcast to $P$ at all, $(ii)$
follows.~$\blacksquare$

\bigskip

We next prove a result similar to Proposition \ref{Prop_Subtree} for
bn-independence in subtrees of a tree.

\begin{theorem}
\label{Prop_subtree_bn}If $T^{\prime}$ is a subtree of a tree $T$, then
$i_{\operatorname{bn}}(T^{\prime})\leq i_{\operatorname{bn}}(T)$.
\end{theorem}

\noindent\textbf{Proof.\hspace{0.1in}}Let $T$ be a tree and $T^{\prime}$ a
subtree of $T$ of order $|V(T)|-k$. Consider an
$i_{\operatorname{bn}}$-broadcast $f$ on $T$. First suppose that $f$ covers
all edges of $T$ and let $P$ be a diametrical path of $T$. By Lemma
\ref{Lem_edge_cover}$(ii)$, 
\begin{equation}
\sigma(f)\geq\left\lceil \frac{1}{2}\ell(P)\right\rceil =\operatorname{rad}%
(T).\label{eq_otherway}%
\end{equation}
Certainly, $\operatorname{rad}(T^{\prime})\leq\operatorname{rad}(T)$ and
$i_{\operatorname{bn}}(T^{\prime})\leq\operatorname{rad}(T^{\prime})$,
therefore, by (\ref{eq_otherway}), 
\[
i_{\operatorname{bn}}(T^{\prime}%
)\leq\mathrm{\operatorname{rad}}(T^{\prime})\leq\mathrm{\operatorname{rad}%
}(T)\leq\sigma(f)=i_{\operatorname{bn}}(T).
\]

Assume therefore that $U_{f}^{E}\neq\varnothing$. Let $\mathcal{C}$ be the
collection of components of $T-U_{f}^{E}$, that is, $\mathcal{C}$ is the
collection of subtrees of $T$ induced by the edges covered by $f$. We
successively delete vertices $v_{1},...,v_{k}$ from $T$ until only $T^{\prime
}$ remains. At each step the vertex that we delete is a leaf of the subtree of $T$ obtained in the previous step, hence each resulting graph is a tree. After each step we adjust the broadcast as necessary so that the
resulting broadcast is maximal bn-independent and has cost no greater than the
cost of $f$. We do this using a recursive algorithm.

Let $T_{0}=T$, $f=f_{0}$ and $\mathcal{C}_{0}=\mathcal{C}$. For $i\geq0$,
assume that we have constructed, in Step $i$, a tree $T_{i}$ that contains
$T^{\prime}$ as subtree, obtained by deleting the vertices $v_{1},...,v_{i}$,
and a maximal bn-independent broadcast $f_{i}$ on $T_{i}$ such that
$\sigma(f_{i})\leq\sigma(f)$. Denote by $\mathcal{C}_{i}$ the collection of
subtrees of $T_{i}$ induced by the edges covered by $f_{i}$.

\noindent\underline{Step $i+1$:}\hspace{0.1in}Let $v_{i}$ be a leaf of
$T_{i}-T^{\prime}$ and $u_{i}$ the stem adjacent to $v_{i}$. Let $H_{i}$ be
the subtree in $\mathcal{C}_{i}$ that contains $v_{i}$. Define $T_{i+1}%
=T_{i}-v_{i}$ and $F_{i+1}=H_{i}-v_{i}$. We consider four cases, depending on
$v_{i}$.\smallskip

\noindent\textbf{Case 1:\hspace{0.1in}}If $f_{i}(v_{i})=0$, let $f_{i+1}$ be
the restriction of $f_{i}$ to $T_{i+1}$. Then $f_{i+1}$ is a bn-independent
and dominating broadcast, and $\sigma(f_{i+1})=\sigma(f_{i})$. The subtree
$H_{i}$ in $\mathcal{C}_{i}$ containing $v_{i}$ becomes the subtree $F_{i+1}$
in $\mathcal{C}_{i+1}$, which has the same number of broadcasting vertices as
$H_{i}$, while all other subtrees in $\mathcal{C}_{i}$ remain unchanged in
$\mathcal{C}_{i+1}$. Since $f_{i}$ is maximal independent, $f_{i+1}$ is
maximal bn-independent, by Proposition \ref{Prop_2bv}.\smallskip

\noindent\textbf{Case 2:\hspace{0.1in}}Suppose $f_{i}(v_{i})\geq2$. Since
$f_{i}$ is bn-independent, $f_{i}(u_{i})=0$. Define the broadcast $f_{i+1}$ on
$T_{i+1}$ by $f_{i+1}(u_{i})=f_{i}(v_{i})-1$ and $f_{i+1}(v)=f_{i}(v)$ for
$v\in V(T_{i+1})-\{u_{i}\}$. Then $\sigma(f_{i+1})<\sigma(f_{i})$ and the
edges of $T_{i+1}$ covered by $f_{i+1}$ are exactly the same as the edges of
$T_{i+1}$ covered by $f_{i}$. Hence $f_{i+1}$ is bn-independent and
dominating. The subtree $F_{i+1}$ in $\mathcal{C}_{i+1}$ has the same number
of broadcasting vertices as $H_{i}$, while all other subtrees remain
unchanged. By Proposition \ref{Prop_2bv}, $f_{i+1}$ is maximal
bn-independent.\smallskip

\noindent\textbf{Case 3:\hspace{0.1in}}Suppose $f_{i}(v_{i})=1$ and $u_{i}$
hears $f_{i}$ from a vertex $w_{i}\neq v_{i}$. Then $w_{i}\in V(H_{i})$.
\label{here}Let $f_{i+1}^{\prime}$ be the restriction of $f_{i}$ to $T_{i+1}$.
Now $\sigma(f_{i+1}^{\prime})=\sigma(f_{i})-1$ and $f_{i+1}^{\prime}$ covers
the same edges of $T_{i+1}$ as $f_{i}$ does, so $f_{i+1}^{\prime}$ is
bn-independent and dominating. However, there are fewer broadcasting vertices
in the component $F_{i+1}$ than in $H_{i}$.

If $F_{i+1}$ has at least two broadcasting vertices, we let $f_{i+1}%
=f_{i+1}^{\prime}$ and proceed as in Case 2. Assume therefore that $w_{i}$ is
the only broadcasting vertex in $F_{i+1}$. If $F_{i+1}=T^{\prime}$, then
$f_{i+1}$ is maximal bn-independent and we are also done. We therefore assume
further that $H_{i}$ is not the only component in $\mathcal{C}_{i}$.

By definition of $F_{i+1}$ and the subtrees in $\mathcal{C}_{i}$, every vertex
in the $f_{i+1}^{\prime}$-boundary of $w_{i}$ also belongs to its private
$f_{i+1}^{\prime}$-boundary. Therefore, the broadcast $f_{i+1}$ on $T_{i+1}$
defined by $f_{i+1}(w_{i})=f_{i}(w_{i})+1$ and $f_{i+1}(v)=f_{i+1}^{\prime
}(v)$ otherwise, is bn-independent, and $\sigma(f_{i+1})=\sigma(f_{i+1}%
^{\prime})+1=\sigma(f_{i})$. Moreover, since $f_{i+1}^{\prime}$ is dominating
and $T_{i+1}$ is connected, there is at least one subtree $G_{i}\in
\mathcal{C}_{i}-\{H_{i}\}$ which contains a vertex that hears the broadcast
$f_{i+1}$ from $w_{i}$. Hence $F_{i+1}$ and $G_{i}$ belong to the same
component in $\mathcal{C}_{i+1}$. Consequently, all subtrees of $T_{i+1}$ in
$\mathcal{C}_{i+1}$ contain at least two broadcasting vertices. By Proposition
\ref{Prop_2bv}, $f_{i+1}$ is maximal bn-independent.\smallskip

\noindent\textbf{Case 4:\hspace{0.1in}}Suppose $f_{i}(v_{i})=1$ and $u_{i}$
belongs to the private $f_{i}$-boundary of $v_{i}$. Then $f_{i}$ does not
cover any edges incident with $u_{i}$ other than $v_{i}u_{i}$. Define
$f_{i+1}$ by $f_{i+1}(u_{i})=1$ and $f_{i+1}(v)=f_{i}(v)$ if $v\neq u_{i}$.
Then $\sigma(f_{i+1})=\sigma(f_{i})$. We now proceed as in Case 3 to show that
$f_{i+1}$ is a maximal bn-independent broadcast on $T_{i+1}$.

\smallskip

It follows that the broadcast $f_{k}$ constructed in Step $k$ is a maximal
independent broadcast on $T_{k}=T^{\prime}$ such that $\sigma(f_{k})\leq
\sigma(f)=i_{\operatorname{bn}}(T)$. Consequently, $i_{\operatorname{bn}%
}(T^{\prime})\leq\sigma(f_{k})\leq i_{\operatorname{bn}}(T)$.~$\blacksquare$

\bigskip

We need another result before we can prove Theorem \ref{Thm_Uniquely_radial}.
We first prove this result for a generalized spider (definition below) and
deduce the required, more general, result by using Theorem
\ref{Prop_subtree_bn}. For $k\geq3$ and $n_{i}\geq1$, $i=1,...,k$, the
\emph{generalized spider} $S=S(n_{1},n_{2},...,n_{k})$ is the tree which has
exactly one vertex $b$ with $\deg(b)=k$, called the \emph{head} of $S$, and
for which the $k$ components of $S-b$ are paths of length $n_{1}%
-1,...,n_{k}-1$, respectively. The paths from $b$ to the end-vertices of $S$
are the \emph{legs} of $S$. The following result by Neilson \cite{N thesis} is
required in the proof of Lemma \ref{Lem_edge_cover_2}.

\begin{theorem}
\label{Thm_paths}\emph{\cite[Theorem 3.1.4]{N thesis} }For any path $P_{n}$
such that $n\neq3$, $i_{\operatorname{bn}}(P_{n})=\left\lceil \frac{2n}%
{5}\right\rceil $.
\end{theorem}

\begin{lemma}
\label{Lem_edge_cover_2}For integers $r\geq2$ and $s\geq1$, consider the
generalized spider $S=S(r-1,r,s)$ with head $b$ and legs $X=(b,x_{1}%
,...,x_{r-1})$, $Y=(b,y_{1},...,y_{r})$ and $Z=(b,z_{1},...,z_{s})$. Suppose
there exists an $i_{\operatorname{bn}}$-broadcast $f$ on $S$ that covers all
edges on $Z$. 

\begin{enumerate}
\item[$(i)$] {Then $\sigma
(f)=i_{\operatorname{bn}}(S)\geq\left\lceil \frac{1}{2}(r+s)\right\rceil $. }
\item[$(ii)$] {If there are $f$-uncovered edges on $X$, then $\sigma
(f)=i_{\operatorname{bn}}(S)\geq\left\lceil \frac{1}{2}(r+s)\right\rceil +1$.}
\end{enumerate}
\end{lemma}

\noindent\textbf{Proof.}\hspace{0.1in}Let $P$ and $Q$ be the paths $X\cup Z$
and $Y\cup Z$, respectively. 
First suppose that $f$ covers all edges on $P$ and hence also on $X$. The bound follows immediately from
Lemma \ref{Lem_edge_cover}$(ii)$ if $f$ overdominates $x_{r-1}$ or $z_{s}$. If
$f$ overdominates neither $x_{r-1}$ nor $z_{s}$, then no broadcasting vertex
on $P$ broadcasts to $y_{r}$. Hence $V_{f}^{+}-V(P)\neq\varnothing$, and again
the result follows from Lemma \ref{Lem_edge_cover}$(ii)$, in which we now have
that $\sigma_{1}+\sigma_{2}\geq1$.

Suppose next that $f$ covers all edges on $Q$. Then $\sigma(f)\geq\left\lceil \frac
{1}{2}\ell(Q)\right\rceil =\left\lceil \frac{1}{2}(r+s)\right\rceil $. If
there is an uncovered edge, say $x_{t}x_{t+1}$, on $X$, then $\{x_{t+2}%
,...,x_{r-1}\}$ contains a vertex in $V_{f}^{+}$; hence $\sum_{v\in
\operatorname{Off}(Q)}f(v)\geq1$. By Lemma \ref{Lem_edge_cover}$(ii)$,
$\sigma(f)\geq\left\lceil \frac{1}{2}\ell(Q)\right\rceil +1=\left\lceil
\frac{1}{2}(r+s)\right\rceil +1$, as required.

\begin{figure}[ht]
    \centering
    \begin{tikzpicture}[main/.style = {draw, circle}, minimum size =0.25 cm, inner sep=0pt, node distance =0.75cm]
    
    \node[main, very thick, fill=yellow] (1) [label={[xshift=0.0cm, yshift=-0.75 cm]1},label={[xshift=-0.55cm, yshift=-0.25cm]\footnotesize $x_{r-1}$}] {};
    \node[rectangle,draw,fill=green] (2) [right of =1] {};
    \node[main,very thick,fill=blue] (3) [right of =2,label={[xshift=0.0cm, yshift=-0.75 cm]\textcolor{blue}{1}}] {};
    \node[main,fill=blue] (4) [right of =3,label={[xshift=0.0cm, yshift=-0.6 cm]\footnotesize $x_{t+1}$}] {};
    \node[main,fill=yellow] (5) [right of =4,label={[xshift=0.0cm, yshift=-0.6 cm]\footnotesize $x_{t}$}] {};
    \node[main,fill=yellow] (6) [right of =5] {};
    \node[main,very thick, fill=yellow] (7) [right of =6,label={[xshift=0.0cm, yshift=-0.75 cm]2},label={[xshift=0.0cm, yshift=0.1 cm]\footnotesize $b$}] {};
    \node[main,fill=yellow] (8) [right of =7] {};
    \node[rectangle,draw,fill=orange] (9) [right of =8] {};
    \node[main,fill=red] (10) [right of =9] {};
    \node[main,fill=red, very thick] (11) [right of =10,label={[xshift=0.0cm, yshift=-0.75 cm]\textcolor{red}{2}}] {};
    \node[main, fill=red] (12) [right of =11] {};
    \node[main,fill=red] (13) [right of =12,label={[xshift=0.4cm, yshift=-0.25cm]\footnotesize $z_s$}] {};
    
    \node[main,fill=yellow] (a) [above of =6,below=0.25cm] {};
    \node[main,fill=yellow] (b) [left of =a, above = 0.2cm,label={[xshift=0.0cm, yshift=0.1 cm]\footnotesize $y_{t'}$}] {};
    \node[main,fill=blue] (c) [left of =b, above = 0.2cm,label={[xshift=0.0cm, yshift=0.1 cm]\footnotesize $y_{t'+1}$}] {};
    \node[very thick, main,fill=blue] (d) [left of =c, above =0.2cm,label={[xshift=0.0cm, yshift=-0.75 cm]\textcolor{blue}{1}}] {};
    \node[rectangle,draw,fill=plum] (e) [left of =d, above = 0.2cm] {};
    \node[very thick,main,fill=red] (f) [left of =e, above =0.2cm,label={[xshift=0.0cm, yshift=-0.75 cm]\textcolor{red}{1}}] {};
    \node[main,fill=red] (g) [left of =f, above = 0.2cm,label={[xshift=-0.4cm, yshift=-0.25cm]\footnotesize $y_{r}$}] {};
    
    \draw(1)--(2);
    \draw(2)--(3);
    \draw(3)--(4);
    \draw[ultra thick](4)--(5);
    \draw(5)--(6);
    \draw(6)--(7);
    \draw(7)--(8);
    \draw(8)--(9);
    \draw(9)--(10);
    \draw(10)--(11);
    \draw(11)--(12);
    \draw(12)--(13);
    \draw(a)--(b);
    \draw[ultra thick](b)--(c);
    \draw(c)--(d);
    \draw(d)--(e);
    \draw(e)--(f);
    \draw(f)--(g);
    \draw(a)--(7);
    
    \draw[<->](-0.1,-1)--(9.2,-1);
    \node[] at (4.55,-1.25) {$P$};
    \draw[->] (4.55,1.25)--(9.2,1.25);
    \draw[->] (4.55,1.25)--(-0.5,3.5);
    \node[] at (4.55,1.7) {$Q$};

    \end{tikzpicture}
    \caption{The spider $S(6,6,7)$ with paths $P$ and $Q$, and uncovered edges $x_tx_{t+1}$ and $y_{t'}y_{t'+1}$}
    \label{Fig_spider}
\end{figure}
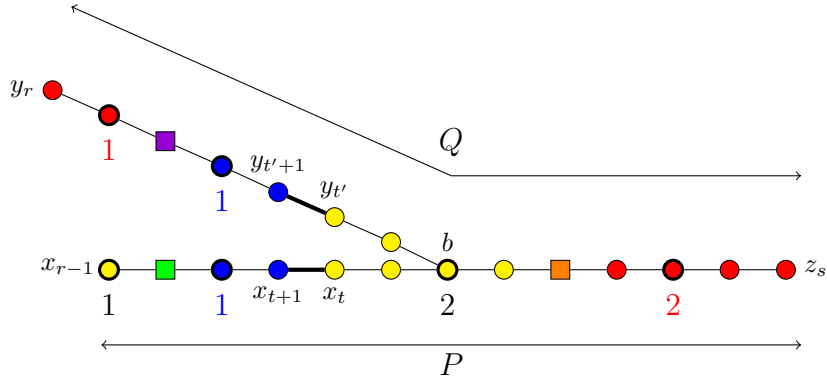

Assume therefore that both $P$ and $Q$ have $f$-uncovered edges. (See Figure
\ref{Fig_spider}.) Since $f$ covers all edges on $Z$, both $X$ and $Y$ have
$f$-uncovered edges. Let $t$ and $t^{\prime}$ respectively be the smallest
indices such that $x_{t}x_{t+1}$ and $y_{t^{\prime}}y_{t^{\prime}+1}$ are
uncovered, and consider the paths
\[
X^{\prime}=(x_{t+1},...,x_{r-1}),\ \ Y^{\prime}=(y_{t^{\prime}+1}%
,...,y_{r}),\ \ R=(x_{t},...,b,z_{1},...,z_{s}),\ \ R^{\prime}=(y_{t^{\prime}%
},...,b,z_{1},...,z_{s}).
\]
Then $|V(X^{\prime})|=r-t-1$ and $|V(Y^{\prime})|=r-t^{\prime}$, while
$\ell(R)=s+t$ and $\ell(R^{\prime})=s+t^{\prime}$. Let $S^{\prime}$ be the
subtree of $S-\{x_{t}x_{t+1},y_{t^{\prime}}y_{t^{\prime}+1}\}$ that contains
$b$ and hence also both paths $R$ and $R^{\prime}$. For each $A\in\{X^{\prime
},Y^{\prime},S^{\prime}\}$, denote the restriction of $f$ to $A$ by $f_{A}$.
By Corollary \ref{Cor_2bv}, $f_{A}$ is a maximal bn-independent broadcast on
$A$, hence $\sigma(f_{A})\geq i_{\operatorname{bn}}(A)$. By definition,
$f_{S^{\prime}}$ covers all edges on $R$ and $R^{\prime}$. Applying Lemma
\ref{Lem_edge_cover}$(ii)$ in turn to $R$ and $R^{\prime}$, we obtain that
\begin{equation}
\sigma(f_{S^{\prime}})\geq\max\left\{  \left\lceil \frac{1}{2}\ell
(R)\right\rceil ,\left\lceil \frac{1}{2}\ell(R^{\prime})\right\rceil \right\}
=\max\left\{  \left\lceil \frac{s+t}{2}\right\rceil ,\left\lceil
\frac{s+t^{\prime}}{2}\right\rceil \right\}  .\label{eq_lengthR}%
\end{equation}
Theorem \ref{Thm_paths} gives%
\begin{equation}
\sigma(f_{X^{\prime}})\geq i_{\operatorname{bn}}(X^{\prime})=\left\lceil
\frac{2(r-t-1)}{5}\right\rceil \text{ and }\sigma(f_{Y^{\prime}})\geq
i_{\operatorname{bn}}(Y^{\prime})=\left\lceil \frac{2(r-t^{\prime})}%
{5}\right\rceil .\label{eq_lengthXY}%
\end{equation}
By Proposition \ref{Prop_2bv}, each component of $T-U_{f}^{E}$ contains at
least two broadcasting vertices. Hence $r-t-1\geq4$ and $r-t^{\prime}\geq4$.
This implies that $\left\lceil \frac{4}{5}(r-t)\right\rceil >\left\lceil
\frac{1}{2}(r-t)\right\rceil $ and $\left\lceil \frac{4}{5}(r-t^{\prime
})\right\rceil >\left\lceil \frac{1}{2}(r-t^{\prime})\right\rceil $. Now, if
$t\geq t^{\prime}$, then $r-t^{\prime}\geq r-t$ and $s+t\geq s+t^{\prime}$.
Substitution in (\ref{eq_lengthR}) and (\ref{eq_lengthXY}) gives%
\begin{align*}
\sigma(f) &  =\sigma(f_{S^{\prime}})+\sigma(f_{X^{\prime}})+\sigma
(f_{Y^{\prime}})\geq\left\lceil \frac{s+t}{2}\right\rceil +\left\lceil
\frac{4(r-t)}{5}\right\rceil \\
&  >\left\lceil \frac{s+t}{2}\right\rceil +\left\lceil \frac{r-t}%
{2}\right\rceil \geq\left\lceil \frac{1}{2}(r+s)\right\rceil .
\end{align*}
On the other hand, if $t<t^{\prime}$, then $s+t^{\prime}>s+t$ and $r-t-1\geq
r-t^{\prime}$, and, from (\ref{eq_lengthR}) and (\ref{eq_lengthXY}),%
\[
\sigma(f)>\left\lceil \frac{s+t^{\prime}}{2}\right\rceil +\left\lceil
\frac{4(r-t^{\prime})}{5}\right\rceil \geq\left\lceil \frac{1}{2}%
(r+s)\right\rceil .
\]
In both cases we see that $i_{\operatorname{bn}}(S)=\sigma(f)\geq\left\lceil \frac{1}{2}%
(r+s)\right\rceil +1$, as required.~$\blacksquare$

\begin{corollary}
\label{Cor_edge_cover2}Let $Q=(u_{0},...u_{s},...,u_{d})$ and $P=(u_{0},...,u_{s},w_{1},$ $...,w_{k})$ be paths of a tree $T$ such that
$u_{0}$ and $w_{k}$ are leaves of $T$, $\{u_{s+1},...,u_{d}\}\cap\{w_{1},$ $...,w_{k}\}=\varnothing$ and $\ell(P)<\ell(Q)$. Let $P^{\prime}$ and $P^{\prime\prime}$ be, respectively,
the $u_{0}-u_{s}$ and $u_{s}-w_{k}$ subpaths of $P$. Suppose there exists an
$i_{\operatorname{bn}}$-broadcast $f$ on $T$ that covers all edges on
$P^{\prime}$.

\begin{enumerate}
\item[$(i)$] Then $\sigma(f)=i_{\operatorname{bn}}(T)\geq\left\lceil \frac{k+s+1}%
{2}\right\rceil >\frac{1}{2}\ell(P)$.

\item[$(ii)$] If some edges on $P^{\prime\prime}$ are $f$-uncovered, then
$\sigma(f)=i_{\operatorname{bn}}(T)\geq\left\lceil \frac{s+k+1}{2}\right\rceil
+1>\frac{1}{2}\ell(P)+1$.
\end{enumerate}
\end{corollary}

\noindent\textbf{Proof.}\hspace{0.1in}Note that $T$ contains the spider
$S=S(k,k+1,s)$ as subtree. Hence, by Theorem \ref{Prop_subtree_bn},
$i_{\operatorname{bn}}(T)\geq i_{\operatorname{bn}}(S)$. Since all conditions
of Lemma \ref{Lem_edge_cover_2} are satisfied for $S$ with respect to its
(maximal) path$\ $of length $k+s$, the result follows from Lemma
\ref{Lem_edge_cover_2}.~$\blacksquare$

\subsection{Outline of proof of Theorem \ref{Thm_Uniquely_radial}}

\label{Sec_outline}Let $T$ be a non-radial tree and $P$ a diametrical path of
$T$. We may think of $T$ as drawn with $P$ on a horizontal line so that we can
refer to a vertex or edge of $P$ as lying to the left (or right) of another
vertex or edge. Assume that the maximum split-set(s) of $T$ has (have)
cardinality $m$, where either $m=1$ or $m=2$. Note that when $m=1$, Definition
\ref{Def_splitset} implies that $T$ has odd diameter, whereas $T$ has even
diameter when $m=2$; that is, $\operatorname{diam}(T)\equiv
m\ (\operatorname{mod}\ 2)$.

In the proof of Theorem \ref{Thm_Uniquely_radial} we consider three types of
vertex partitions of $T$ into subtrees. First, for a maximum split-set $M$ of
$T$, there is the associated partition of $T$ into $m+1$ uniquely radial
subtrees. Then, in order to obtain a contradiction, we assume that
$i_{\operatorname{bn}}(T)=\gamma_{b}(T)$ instead of $\gamma_{b}(T)+1$, and
consider a maximal bn-independent broadcast $f$ of cost $\sigma
(f)=i_{\operatorname{bn}}(T)=\gamma_{b}(T)$. We show that $f$ leaves some
edges $e_{1},...,e_{k}$ of $P$ uncovered, and consider the subtrees
$T_{1},...,T_{k+1}$ of $T$ obtained by deleting $e_{1},...,e_{k}$, as well as
the restrictions $f_{i}$ of $f$ to $T_{i}$. Then we consider a succession of
partitions of $T$ into two subtrees -- the trees $L_{i}$ to the left of
$e_{i}$, and $R_{i}$ to the right, for $i=1,...,k$, together with the
restrictions $g_{i}$ and $h_{i}$ of $f$ to $L_{i}$ and $R_{i}$, respectively.
As we proceed from left to right, we iteratively determine $\sigma(g_{i})$.
Our aim is to construct a strictly increasing sequence of indices $i_{1}%
,i_{2},...$ of infinite length such that $\sigma(g_{i_{j}})\geq
\operatorname{diam}(L_{i_{j}})$, thus obtaining a contradiction. To this end
we also consider a variety of subtrees of $L_{i}$ and $R_{i}$.

We state and prove a number of claims within the proof; the end of the proof
of each claim is indicated by an open diamond ($\lozenge$). When we consider
different cases, the end of the proof of each case is indicated by a solid
diamond ($\blacklozenge$).

\subsection{Proof of Theorem \ref{Thm_Uniquely_radial}}

\label{Sec_unique}We restate the theorem for convenience.

\bigskip

\noindent\textbf{Theorem \ref{Thm_Uniquely_radial}}\hspace{0.1in}\emph{Let
}$T$\emph{ be a tree whose maximum split-set(s) have cardinality }%
$m\in\{1,2\}$\emph{. Suppose that for each maximum split-set }$M$\emph{, every
component of }$T-M$\emph{ is uniquely radial. Then }$i_{\operatorname{bn}%
}(T)=\gamma_{b}(T)+1$\emph{.}

\bigskip

\noindent\textbf{Proof}.\hspace{0.1in}Suppose, contrary to the statement, that
$i_{\operatorname{bn}}(T)=\gamma_{b}(T)$. Let $P=(v_{0},v_{1},...,v_{d})$ be a
diametrical path of $T$ and consider an $i_{\operatorname{bn}}$-broadcast $f$
on $T$. By Theorem \ref{Thm_gammaT},
\begin{equation}
\sigma(f)=\gamma_{b}(T)=\frac{1}{2}(\operatorname{diam}(T)-m)<\frac{1}%
{2}\operatorname{diam}(T). \label{eq_gamma_b}%
\end{equation}
By Lemma \ref{Lem_edge_cover}, $f$ does not cover all edges of $P$. Let
$e_{1},...,e_{k}$, $k\geq1$, be the uncovered edges on $P$; say $e_{i}%
=x_{i}y_{i}$ (where $x_{i}$ is to the left of $y_{i}$), and also let
$y_{0}=v_{0}$ and $x_{k+1}=v_{d}$. Note that is is possible that
$x_{i+1}=y_{i}$ for some $i$; this happens precisely when some vertex in
$\operatorname{Touch}(P)-V(P)$ broadcasts to exactly one vertex on $P$. Let
$T_{1},...,T_{k+1}$ be the components of $T-\{e_{1},...,e_{k}\}$, where
$T_{i}$ is the component whose leftmost and rightmost vertices on $P$ are
$y_{i-1}$ and $x_{i}$, respectively. For each $i$, let $P_{i}$ be the subpath
of $P$ on $T_{i}$, and $f_{i}$ the restriction of $f$ to $T_{i}$. Since
$k\geq1$, $|V_{f}^{+}|\geq2$. By Corollary \ref{Cor_2bv}, each $f_{i}$ is a
maximal bn-independent broadcast on $T_{i}$, hence $\sigma(f_{i})\geq
i_{\operatorname{bn}}(f_{i})$ for each $i$. By Proposition \ref{Prop_2bv},
$|V_{f_{i}}^{+}|\geq2$ for each $i$. 

Suppose $f_{i}$ is not an $i_{\operatorname{bn}}$-broadcast and let $g_{i}$ be
an $i_{\operatorname{bn}}$-broadcast on $T_{i}$ instead. Then $\sigma
(g_{i})<\sigma(f_{i})$. By maximality, $g_{i}$ is a dominating broadcast on
$T_{i}$ (Proposition~\ref{Thm_Minimal}). However, this implies that
$g=(f-f_{i})\cup g_{i}$ is a dominating broadcast on $T$ with cost
$\sigma(g)<\sigma(f)=\gamma_{b}(T)$, a contradiction. Hence $f_{i}$ is an
$i_{\operatorname{bn}}$-broadcast on $T_{i}$ for each $i$.

For $i=1,...,k$, let $L_{i}$ and $R_{i}$ ($L$ for \textquotedblleft
left\textquotedblright, $R$ for \textquotedblleft right\textquotedblright) be
the subtrees of $T-e_{i}$ that contain $x_{i}$ and $y_{i}$, respectively. Let
$P(L_{i})$ and $P(R_{i})$ be the subpaths of $P$ on $L_{i}$ and $R_{i}$,
respectively. Also denote the restriction of $f$ to $L_{i}$ and $R_{i}$ by
$g_{i}$ and $h_{i}$, respectively. As in the case of the $f_{i}$, the $g_{i}$
and $h_{i}$ are $i_{\operatorname{bn}}$-broadcasts on $L_{i}$ and $R_{i}$, respectively.

\begin{claim}
\label{Cl_1}None of the uncovered edges $e_{1},...,e_{k}$ is a split-edge of
$T$.
\end{claim}

\noindent\textbf{Proof of Claim \ref{Cl_1}.\hspace{0.1in}}Suppose that for
some $j$, $j=1,...,k$, $e_{j}$ is a split-edge of $T$ on $P$. Consider $L_{j}$
and $R_{j}$, and note that $g_{j}=\bigcup_{i=1}^{j}f_{i}$ and $h_{j}%
=\bigcup_{i=j+1}^{k+1}f_{i}$. Since $|V_{f_{i}}^{+}|\geq2$ for each $i$,
$|V_{g_{j}}^{+}|,|V_{h_{j}}^{+}|\geq2$. Now, if $m=1$, then for $M=\{e_{j}\}$,
both $L_{j}$ and $R_{j}$ are radial subtrees of $T$. On the other hand, if
$m=2$, then $M=\{e_{j},e^{\prime}\}$ for some $e^{\prime}\in E(P)$. We may
assume without loss of generality that $e_{j}$ lies to the left of $e^{\prime
}$ on $P$. Then $L_{j}$ is a radial subtree of $T$. In either case, $L_{j}$ is
uniquely radial, which contradicts $|V_{g_{j}}^{+}|\geq2$. We conclude that no
$e_{i}$ is a split-edge.$~\lozenge$

\bigskip

We continue with the proof of Theorem \ref{Thm_Uniquely_radial} and consider
two cases, depending on whether $k=1$ or $k\geq2$.

\bigskip

\noindent\textbf{Case 1:\hspace{0.1in}}$k=1$. Then $\operatorname{diam}%
(T)=\ell(P_{1})+\ell(P_{2})+1$. As stated in (\ref{eq_gamma_b}), $\gamma
_{b}(T)=\frac{1}{2}(\operatorname{diam}(T)-m)$. On the other hand, $\gamma
_{b}(T)=\sigma(f_{1})+\sigma(f_{2})$. Since $f_{i}$ covers each edge of
$P_{i}$, Lemma \ref{Lem_edge_cover}$(ii)$ implies that $\sigma(f_{i})\geq
\ell(P_{i})/2$ for $i=1,2$. Therefore%
\begin{equation}
\operatorname{diam}(T)-m=2(\sigma(f_{1})+\sigma(f_{2}))\geq\ell(P_{1}%
)+\ell(P_{2})=\operatorname{diam}(T)-1. \label{eq_diam1}%
\end{equation}
We conclude that $m=1$ and equality holds throughout (\ref{eq_diam1}). In
particular,
\begin{equation}
\gamma_{b}(T_{i})=\frac{1}{2}\ell(P_{i})\text{\ for\ }i=1,2, \label{eq_k=1a}%
\end{equation}
hence $\ell(P_{i})$ is even. Since $e_{1}$ is not a split-edge, Definition
\ref{Def_splitset} implies that $P_{i}$ is not a diametrical path of $T_{i}$
for at least one $i\in\{1,2\}$. Assume without loss of generality that $P_{1}$
is not a diametrical path of $T_{1}$. Since $f_{1}$ is an
$i_{\operatorname{bn}}$-broadcast on $T_{1}$ that covers all edges on $P_{1}$,
Corollary \ref{Cor_edge_cover2} implies that $i_{\operatorname{bn}}%
(T_{1})>\frac{1}{2}\ell(P_{1})$, which contradicts (\ref{eq_k=1a}). This
completes the proof of Case~1.$~\blacklozenge$

\bigskip

\noindent\textbf{Case 2:\hspace{0.1in}}$k\geq2$. We first state and prove
another claim.

\begin{claim}
\label{Cl_2}Each $P_{i},\ i=2,...,k$, has even length.
\end{claim}

\noindent\textbf{Proof of Claim \ref{Cl_2}.\hspace{0.1in}}Since the edges
$e_{i}$ are uncovered, no $x_{i}$ or $y_{i}$ is overdominated, for
$i=1,...,k$. Since $f$ is bn-independent, each edge of $P_{i}$ is covered
exactly once. Since each vertex in $V_{f_{i}}^{+}$ covers an even number of
edges on $P_{i}$, the result follows.~$\lozenge$

\bigskip

The same argument shows that if $P_{1}$ or $P_{k+1}$ has odd length, then
$y_{0}$ or $x_{k+1}$, respectively, is overdominated. Since every edge of
$P_{i}$ is covered, Lemma \ref{Lem_edge_cover}$(ii)$ implies that, for each
$i$,
\begin{equation}
\sigma(f_{i})\geq\left\lceil \frac{1}{2}\ell(P_{i})\right\rceil .
\label{eq_sigma_Fi}%
\end{equation}
Since
\[
\sum_{i=1}^{k+1}\sigma(f_{i})=\gamma_{b}(T)=\frac{1}{2}(\operatorname{diam}%
(T)-m)=\frac{1}{2}\left[  \sum_{i=1}^{k+1}\ell(P_{i})+k-m\right]  ,
\]
we see that%
\[
2\sum_{i=1}^{k+1}\sigma(f_{i})-\sum_{i=1}^{k+1}\ell(P_{i})=k-m\geq0.
\]
Hence either $k=m=2$ and $\sigma(f_{i})=\frac{1}{2}\ell(P_{i})$ for each $i$,
or $k-m>0$ and there is at least one value of $i$ for which $\sigma
(f_{i})>\frac{1}{2}\ell(P_{i})$. In the former case we proceed as in Case 1
(but with $m=2$) to obtain a contradiction. Hence we may assume that the
latter case holds. For such an $i$, since each edge of $P_{i}$ is covered
exactly once, either (a) $V_{f_{i}}^{+}-V(P_{i})\neq\varnothing$, or (b)
$i\in\{1,k+1\}$, $\ell(P_{i})$ is odd and an end-vertex of $P_{i}$ is
overdominated (or both (a) and (b) hold).

We consider two subcases, depending on whether $\ell(P_{1})$ is odd or even,
beginning with the case where $\ell(P_{1})$ is odd. The analysis for this case
can be repeated with only a trivial modification when $\ell(P_{1})$ is even.

\bigskip

\noindent\textbf{Case 2(a):\hspace{0.1in}}Suppose $\ell(P_{1})$ is odd. By
(\ref{eq_sigma_Fi}), $\sigma(f_{1})\geq\frac{1}{2}(\ell(P_{1})+1)$. Since
$\ell(P_{1})$ is odd while $\ell(P_{2})$ is even, and $P(L_{2})$ consists of
$P_{1}$ followed by $e_{1}$ followed by $P_{2}$, $P(L_{2})$ has even length
$\ell(P_{1})+\ell(P_{2})+1$. Moreover,
\begin{equation}
\sigma(g_{2})=\sigma(f_{1})+\sigma(f_{2})\geq\frac{1}{2}[\ell(P_{1}%
)+\ell(P_{2})+1]=\frac{1}{2}\ell(P(L_{2})). \label{eq_g2}%
\end{equation}

\begin{figure}[ht]
    \centering
    \begin{tikzpicture}[main/.style = {draw, circle}, minimum size =0.25 cm, inner sep=0pt, node distance =0.75cm]
    
    \node[main,fill=blue,very thick] (1) [label={[xshift=-0.4cm, yshift=-0.25cm]\footnotesize $v_0$}] {};
    \node[rectangle,draw,fill=green] (2) [right of =1] {};
    \node[main,fill=yellow] (3) [right of =2]{};
    \node[main,fill=yellow,very thick] (4) [right of =3]{};
    \node[main, fill=yellow] (5) [right of =4]{};
    \node[main,fill=yellow] (6) [right of =5,label={[xshift=0cm, yshift=0.025cm]\footnotesize $x_1$}]{};
    \node[main,fill=blue] (7) [right of =6,label={[xshift=0cm, yshift=0.025cm]\footnotesize $y_2$}] {};
    \node[fill=blue,main] (8) [right of =7,very thick]{};
    \node[draw,rectangle,fill=plum] (9) [right of =8] {};
    \node[main,fill=red,very thick] (10) [right of =9,label={[xshift=0.2cm, yshift=-0.03cm]\footnotesize $z_2$}] {};
    \node[main,fill=red] (11) [right of =10,label={[xshift=0cm, yshift=0.025cm]\footnotesize $x_2$}] {};
    \node[main] (12) [right of =11,label={[xshift=0cm, yshift=0.025cm]\footnotesize $y_2$}]{};
    \node[main,fill=black] (13) [right of =12, right=0.5cm,label={[xshift=0.4cm, yshift=-0.25cm]\footnotesize $v_d$}] {};
    
    \node[draw,rectangle,fill=orange] (14) [above of =5] {};
    \node[main,fill=red,very thick] (15) [above of =14] {};
    \node[main,fill=red] (16) [above of =15] {};
    \node[main,fill=red] (17) [above of =4, above = 1cm] {};
    \node[fill=blue,main,above of =8] (18) [left=0.25cm]{};
    \node[fill=blue,main, above of =8] (19) [right=0.25cm] {};
    
    \node[rectangle,fill=orange,draw] (20) [above of =10] {};
    \node[main,very thick,fill=yellow] (21) [above of =20] {};
    \node[main,fill=yellow] (22) [above of =21,label={[xshift=0cm, yshift=0.1cm]\footnotesize $q_2$}]{};
    
    \draw(1)--(2);
    \draw(2)--(3);
    \draw(3)--(4);
    \draw(4)--(5);
    \draw(5)--(6);
    \draw[ultra thick](6)--(7);
    \draw(7)--(8);
    \draw(8)--(9);
    \draw(9)--(10);
    \draw(10)--(11);
    \draw[ultra thick](11)--(12);
    \draw[dashed] (12)--(13);
    \draw(14)--(5);
    \draw(15)--(14);
    \draw(15)--(17);
    \draw(15)--(16);
    \draw(8)--(18);
    \draw(8)--(19);
    \draw(10)--(20);
    \draw(20)--(21);
    \draw(21)--(22);
    
    \draw[] (-0.65,-0.5) rectangle (7.7,3);
    \draw[] (8.05,-0.5) rectangle (10.25,3);
     {[on background layer]\draw[fill=orange,opacity=0.2] (-0.2,-0.4) rectangle (4,2.8);
    \draw[fill=orange,opacity=0.2] (4.25,-0.4) rectangle (7.65,2.8);}
    
    \draw[<->](-0.2,-0.75)--(7.65,-0.75);
    \node[] at (3.725,-1.15) {$P(L_2)$};
    \node[] at (0.2,2.4) {$T_1$};
    \node[] at (4.7 ,2.4) {$T_2$};
    \node[] at (7.25,3.3) {$L_2$};
    \node[] at (9.7,3.3 ) {$R_2$};
    \draw[<-] (-0.1,0.45)--(6.25,0.45);
    \draw[->] (6.25,0.45)--(6.25,2.4);
    \node[] at (1,0.75) {$Q(L_2)$};
    
    \end{tikzpicture}
    \caption{A step in the proof of Theorem 2.6}
    \label{Fig_eq_bound}
\end{figure}
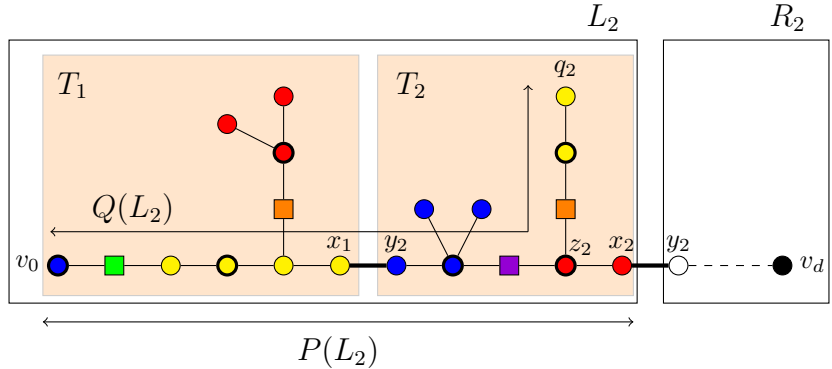

Our next goal is to show that there exists an even integer, say $2t$, such
that $P(L_{2t})$ is a diametrical path of $L_{2t}$. (See Figure
\ref{Fig_eq_bound}.)

\begin{claim}
\label{Cl_3}For\smallskip\ $p\in\{1,...,\left\lfloor \frac{k}{2}\right\rfloor
\}$, suppose $P(L_{2i})$ is not a diametrical path of $L_{2i}$ for each
$i=1,...,p-1$.\smallskip\ Then $\sigma(g_{2p})\geq\frac{1}{2}\ell(P(L_{2p}))$.
If, in addition, $P(L_{2p})$ is not a diametrical path of $L_{2p}$, then
$\sigma(g_{2p})\geq\frac{1}{2}\ell(P(L_{2p}))+1$.
\end{claim}

\noindent\textbf{Proof of Claim \ref{Cl_3}.\hspace{0.1in}}Let $Q(L_{2i})$ be a
diametrical path of $L_{2i},\ i=1,...,p$. Since $v_{0}=y_{0}$ is a peripheral
vertex of $T$, we may assume without loss of generality that $v_{0}$ is also a
peripheral vertex of $Q(L_{2i})$. Let $z_{2i}$ be the vertex nearest to
$x_{2i}$ common to $Q(L_{2i})$ and $P(L_{2i})$, and let $q_{2i}$ be the
end-vertex of $Q(L_{2i})$ antipodal to $v_{0}$. By assumption, $d(z_{2i}%
,q_{2i})>d(z_{2i},x_{2i})$. If $z_{2i}\in V(T_{1})$, then $g_{2i}$ covers
every edge on the $y_{0}-z_{2i}$ subpath of $P(L_{2i})$. By Corollary
\ref{Cor_edge_cover2}, $i_{\operatorname{bn}}(L_{2i})=\sigma(g_{2i})>\frac
{1}{2}\ell(P(L_{2i}))$, and since $\frac{1}{2}\ell(P(L_{2i}))$ is even,
$\sigma(g_{2i})\geq\frac{1}{2}\ell(P(L_{2i}))+1$. We therefore assume
henceforth that $z_{2i}\notin V(T_{1})$.

Suppose $p=1$. We know from (\ref{eq_g2}) that $\sigma(g_{2})\geq\frac{1}%
{2}\ell(P(L_{2}))$, hence the first part of the claim holds. Assume
$\ell(Q(L_{2}))>\ell(P(L_{2}))$. Since $z_{2}\notin V(T_{1})$, we deduce that
$z_{2}\in V(T_{2})$. Since $f_{2}$ covers every edge on the $y_{1}-z_{2}$
subpath of $P_{2}$, Corollary \ref{Cor_edge_cover2} implies that
$i_{\operatorname{bn}}(T_{2})=\sigma(f_{2})>\frac{1}{2}\ell(P_{2})$. Since
$\ell(P_{2})$ is even, $\sigma(f_{2})\geq\frac{1}{2}\ell(P_{2})+1$. Therefore
$\sigma(g_{2})=\sigma(f_{1})+\sigma(f_{2})\geq\frac{1}{2}(\ell(P_{1}%
)+\ell(P_{2})+1)+1$, that is,
\begin{equation}
\sigma(g_{2})\geq\frac{1}{2}\ell(P(L_{2}))+1, \label{eq_L2a}%
\end{equation}
and the second part of the claim holds.

Suppose $p\geq2$. Assume that we have shown the result for $1,...,p-1$, and
consider $L_{2p}$. Since $\ell(P(L_{2p}))=\ell(P(L_{2p-2}))+\ell
(P_{2p-1})+\ell(P_{2p})+2$, our assumption and (\ref{eq_sigma_Fi}) prove the
first part of the claim, namely%
\begin{equation}
\sigma(g_{2p})=\sigma(g_{2p-2})+\sigma(f_{2p-1})+\sigma(f_{2p})\geq\frac{1}%
{2}[\ell(P(L_{2p-2})+2+\ell(P_{2p-1})+\ell(P_{2p})]=\frac{1}{2}\ell
(P(L_{2p})). \label{eq_L2p}%
\end{equation}

Now assume $\ell(Q(L_{2p})>\ell(P(L_{2p}))$. Then $z_{2p}\in V(T_{\alpha})$
for some $\alpha\in\{2,3,...,2p\}$. Let $G$ be the subtree of $L_{2p}%
-x_{\alpha-1}y_{\alpha-1}$ that contains $y_{\alpha-1}$ and hence also
$z_{2p}$. Let $P(G)$ be the subpath of $P$ in $G$, and $g$ the restriction of
$f$ to $G$. As in the case of the $f_{i}$ and $g_{i}$, the broadcast $g$ is
maximal bn-independent. We may apply Corollary \ref{Cor_edge_cover2} to $G$,
with the $y_{\alpha-1}-z_{2p}$ subpath of $P(L_{2p})$ as $P^{\prime}$, to
obtain that $\sigma(g)>\frac{1}{2}\ell(P(G))$.

\begin{itemize}
\item Suppose $\alpha=2$. Then $\ell(P(L_{2p}))=\ell(P_{1})+\ell(P(G))+1$ and
$\ell(P(G))$ is even (since $\ell(P(L_{2p}))$ is even and $\ell(P_{1})$ is
odd), hence $\sigma(g)\geq\frac{1}{2}\ell(P(G))+1$. Therefore
\[
\sigma(g_{2p})=\sigma(f_{1})+\sigma(g)\geq\frac{1}{2}\left[  \ell
(P_{1})+1\right]  +\frac{1}{2}[\ell(P(G))+2]=\frac{1}{2}\ell(P(L_{2p}))+1.
\]

\item Suppose $\alpha>2$ is even. Then
\[
\ell(P(L_{2p}))=\ell(P(L_{\alpha-2}))+\ell(P_{\alpha-1})+\ell(P(G))+2.
\]
It follows that $\ell(P(G))$ is even, hence $\sigma(g)\geq\frac{1}{2}%
\ell(P(G))+1$. By assumption, $\sigma(g_{\alpha-2})\geq\frac{1}{2}%
\ell(P(L_{\alpha-2}))+1$, hence
\[
\sigma(g_{2p})\geq\frac{1}{2}\left[  \ell(P(L_{\alpha-2}))+\ell(P_{\alpha
-1})+\ell(P(G))\right]  +2=\frac{1}{2}\ell(P(L_{2p}))+1.
\]

\item Suppose $\alpha\geq3$ is odd. Then
\[
\ell(P(L_{2p}))=\ell(P(L_{\alpha-1})+\ell(P(G))+1
\]
and $\ell(P(G))$ is odd, so $\sigma(g)\geq\frac{1}{2}[\ell(P(G))+1]$. By
assumption, $\sigma(g_{\alpha-1})\geq\frac{1}{2}\ell(P(L_{\alpha-1}))+1$.
Therefore%
\begin{align*}
\sigma(g_{2p})  &  =\sigma(g_{\alpha-1})+\sigma(g)\\
&  \geq\frac{1}{2}[\ell(P(L_{\alpha-1}))+\ell(P(G))+1]+1\geq\frac{1}{2}%
\ell(P(L_{2p}))+1.~~\lozenge
\end{align*}

\end{itemize}

\bigskip

Applying Claim \ref{Cl_3} to $p=\left\lfloor \frac{k}{2}\right\rfloor $, we
obtain that $\sigma(g_{2p})\geq\frac{1}{2}\ell(P(L_{2p}))+1$. We show next
that this result produces a contradiction.

\begin{claim}
\label{Cl_4}Suppose that $P(L_{2i})$ is not a diametrical path of $L_{2i}$ for
each $i=1,...,\left\lfloor \frac{k}{2}\right\rfloor $. Then $\gamma_{b}%
(T)\geq\frac{1}{2}\operatorname{diam}(T)$.
\end{claim}

\noindent\textbf{Proof of Claim \ref{Cl_4}.\hspace{0.1in}}Suppose $k$ is even.
Then $2p=k$ and $L_{2p}$ is followed by $T_{k+1}=T_{2p+1}$, hence
\[
\operatorname{diam}(T)=\ell(P)=\ell(P(L_{2p})+\ell(P_{2p+1})+1.
\]
By (\ref{eq_sigma_Fi}) and Claim \ref{Cl_3},
\[
\gamma_{b}(T)=\sigma(g_{2p})+\sigma(f_{2p+1})\geq\frac{1}{2}[\ell
(P(L_{2p}))+\ell(P_{2p+1})]+1>\frac{1}{2}\operatorname{diam}(T).
\]

Suppose $k$ is odd. Then $2p=k-1$ and $L_{2p}$ is followed by $T_{k}=T_{2p+1}$
and $T_{k+1}=T_{2p+2}$, hence
\[
\operatorname{diam}(T)=\ell(P)=\ell(P(L_{2p})+\ell(P_{2p+1})+\ell
(P_{2p+2})+2,
\]
and%
\begin{align*}
\gamma_{b}(T)  &  =\sigma(g_{2p})+\sigma(f_{2p+1})+\sigma(f_{2p+2})\\
&  \geq\frac{1}{2}[\ell(P(L_{2p}))+\ell(P_{2p+1})+\ell(P_{2p+2})+2]=\frac
{1}{2}\operatorname{diam}(T).~~\lozenge
\end{align*}

\bigskip

However, Claim \ref{Cl_4} contradicts (\ref{eq_gamma_b}). We conclude that
there exists a smallest even integer $2t$, $t\geq1$, such that $P(L_{2t})$ is
a diametrical path of $L_{2t}$. By Claim \ref{Cl_3} and the choice of $t$,
$\sigma(g_{2t})\geq\frac{1}{2}\ell(P(L_{2t}))=\frac{1}{2}\operatorname{diam}%
(L_{2t})$. Since $\sigma(f)<\frac{1}{2}\operatorname{diam}(T)$, the edge
$e_{2t}=x_{2t}y_{2t}$ and the tree $R_{2t}$ exist.

By Claim \ref{Cl_1}, $e_{2t}$ is not a split-edge of $T$; that is, there is no
split-set $M$ such that $M=\{e_{2t}\}$ (if $m=1)$ or $M=\{e_{2t},e^{\prime}\}$
for some $e^{\prime}\in E(P)$ (if $m=2$). Hence Definition \ref{Def_splitset}
and the fact that $\ell(P(L_{2t})$ is even imply that either

\begin{enumerate}
\item[(a)] $P(R_{2t})=y_{2t},...,v_{d}$ is not a diametrical path of $R_{2t}$, or

\item[(b)] $m=2$ and $P(R_{2t})$ is a diametrical path of $R_{2t}$, but there
is no edge $e^{\prime}$ on $P(R_{2t})$ such that $\{e_{2t},e^{\prime}\}$ is a
split-set of $T$.
\end{enumerate}

Suppose (b) holds. Then $\operatorname{diam}(T)$ is even, hence
$\operatorname{diam}(R_{2t})$ is odd. Since any split-edge $e^{\prime\prime}$
of $R_{2t}$ on its diametrical path $P(R_{2t})$ would result in $\{e_{2t}%
,e^{\prime\prime}\}$ being a split-set of $T$, which is not the case, we
conclude that $R_{2t}$ is a bicentral radial tree. But then $\gamma_{b}%
(R_{2t})=\frac{1}{2}(\ell(P(R_{2t}))+1)$, thus
\[
\gamma_{b}(T)=\gamma_{b}(L_{2t})+\gamma_{b}(R_{2t})\geq\frac{1}{2}%
[\ell(P(L_{2t}))+(\ell(P(R_{2t}))+1)]>\frac{1}{2}(\operatorname{diam}(T)-m),
\]
a contradiction. Therefore (a) holds. We proceed to show that there is a
smallest integer $t_{1}>t$ such that $P(L_{2t_{1}})$ is a diametrical path of
$L_{2t_{1}}$ and $\sigma(g_{2t_{1}})\geq\frac{1}{2}\ell(P(L_{2t_{1}}))$. This
proof first considers the subtrees $R_{i}$ instead of $L_{i}$ and shows some
similarities with the proof of Claim \ref{Cl_3}.

\begin{figure}[ht]
    \centering
    \begin{tikzpicture}[main/.style = {draw, circle}, minimum size =0.25 cm, inner sep=0pt, node distance =0.75cm]
    
    \node[main] (1) [label={[xshift=0cm, yshift=-0.6cm]\footnotesize $y_0$}] {};
    \node[main] (2) [right of =1,right= 1cm,label={[xshift=0cm, yshift=-0.6cm]\footnotesize $x_{2t}$}] {};
    \node[main] (3) [right of =2,label={[xshift=0.1cm, yshift=-0.6cm]\footnotesize $y_{2t}$}] {};
    \node[main] (4) [right of =3, right = 1cm] {};
    \node[main,fill=plum] (5) [right of =4,label={[xshift=0.4cm, yshift=-0.6cm]\footnotesize $y_{2t+\alpha-1}$}] {};
    \node[main,fill=plum] (6) [right of =5,right=1cm,label={[xshift=0.4cm, yshift=0.1cm]\footnotesize $w_{2t}$}] {};
    \node[main,fill=plum] (7) [right of =6,right=1cm,label={[xshift=-0.3cm, yshift=-0.6cm]\footnotesize $x_{2t+\alpha}$}] {};
    \node[main] (8) [right of =7] {};
    \node[main] (9) [right of =8, right = 1cm,label={[xshift=0cm, yshift=-0.6cm]\footnotesize $v_{d}$}] {};
    \node[main,fill=plum] (10) [above of =6, above=1cm,label={[xshift=0cm, yshift=0.1cm]\footnotesize $u_{2t}$}]{};
    
    \draw[dashed] (1)--(2);
    \draw (2) -- (3);
    \draw[dashed] (3)--(4);
    \draw(4)--(5);
    \draw[dashed](5)--(6);
    \draw[dashed] (6)--(7);
    \draw(7)--(8);
    \draw[dashed](8)--(9);
    \draw[dashed] (10)--(6);
    
    \node[] at (-0.5,0.1) {$P$};
    \node[] at (5.75,2) {\textcolor{plum}{$T_{2t+\alpha}$}};
    \node[] at (5.85,-1.05) {\textcolor{red}{$P(H)$}};
    \node[] at (2.8,2.25) {\textcolor{red}{$H$}};
    \node[] at (0,3.25) {$L_{2t+\alpha}$};
    \node[] at (11.45,2.8) {\textcolor{blue}{$R_{2t}$}};
    \node[] at (2.2,0.2) {\footnotesize{$e_{2t}$}};
    \node[] at (5,0.3) {\footnotesize{$e_{2t+\alpha-1}$}};
    \node[] at (9.7,0.3) {\footnotesize{$e_{2t+\alpha}$}};
    \node[] at (10.1,1.3) {\textcolor{blue}{$Q(R_{2t})$}};

    \draw[] (-0.8,-2.5) rectangle (9.25,3.75); 
    \draw[red] (2.45,-1.5) rectangle (9.25,2.75);
    \draw[<->,red] (2.6,-1.3)--(9.1,-1.3);
    \draw[blue] (2.45, -2) rectangle (12,3.25);
    {[on background layer]\draw[fill=orange, opacity=0.2] (5,-0.75) rectangle (9.25, 2.5);}
    \draw[->,blue] (8,1) -- (11.75,1);
    \draw[->,blue] (8,1) -- (8,2);

    \end{tikzpicture}
    \caption{The subtree $R_{2t}$, its diametrical path $Q(R_{2t})$, and the subtree $H$ of $R_{2t}$}
    \label{Fig_MainH}
\end{figure}

The next part of the proof is illustrated in Figure~\ref{Fig_MainH}. Let
$Q(R_{2t})$ be a diametrical path of $R_{2t}$. Since $v_{d}=x_{k+1}$ is a
peripheral vertex of $T$, we may assume without loss of generality that
$v_{d}$ is a peripheral vertex of $Q(R_{2t})$ as well. Let $u_{2t}$ be the
end-vertex of $Q(R_{2t})$ antipodal to $v_{d}$ and let $w_{2t}$ be the vertex
nearest to $y_{2t}$ that is common to $Q(R_{2t})$ and $P(R_{2t})$. Say
$w_{2t}\in V(T_{2t+\alpha})$ for some integer $\alpha\geq1$. If $2t+\alpha
=k+1$, let $H=R_{2t}$. Otherwise, let $H$ be the subtree of $R_{2t}%
-e_{2t+\alpha}$ that contains $x_{2t+\alpha}$. Denote the subpath of $P$ in
$H$ by $P(H)$, and the restriction of $f$ to $H$ by $h$. By Corollary
\ref{Cor_2bv}, $h$ is a maximal bn-independent broadcast on $H$. Moreover, $h$
covers every edge on the $w_{2t}-x_{2t+\alpha}$ subpath of $P(H)$. By
Corollary \ref{Cor_edge_cover2}, $\sigma(h)>\frac{1}{2}\ell(P(H))$.

It follows that $R_{2t}\neq H$, otherwise we have $\operatorname{diam}%
(T)=\ell(P)=\ell(P(L_{2t}))+\ell(P(H))+1$ and $\sigma(f)=\sigma(g_{2t}%
)+\sigma(h)>\frac{1}{2}[\ell(P(L_{2t}))+\ell(P(H))]\geq\frac{1}{2}%
(\operatorname{diam}(T)-m)$, which is not the case. Hence we can consider the
trees $L_{2t+\alpha}$ and $R_{2t+\alpha}$. Note that $\ell(P(L_{2t+\alpha
}))=\ell(P(L_{2t}))+\ell(P(H))+1$.

\begin{itemize}
\item Suppose $\alpha>1$. Then the edge $x_{2t+\alpha-1}y_{2t+\alpha-1}$ on
the $w_{2t}-y_{2t}$ subpath of $P(H)$ is uncovered. By the second part of
Corollary \ref{Cor_edge_cover2}, $\sigma(h)>\frac{1}{2}\ell(P(H))+1$. If
$\ell(P(L_{2t+\alpha}))$ is even, i.e., if $\alpha$ is even, then $\ell(P(H))$
is odd and $\sigma(h)\geq\frac{1}{2}(\ell(P(H))+3)$. Hence
\[
\sigma(g_{2t+\alpha})=\sigma(g_{2t})+\sigma(h)\geq\frac{1}{2}[\ell
(P(L_{2t}))+\ell(P(H))+3]=\frac{1}{2}\ell(P(L_{2t+\alpha}))+1.
\]
Thus, whether or not $P(L_{2t+\alpha})$ is a diametrical path of
$L_{2t+\alpha}$, we have that $\sigma(g_{2t+\alpha})\geq\frac{1}{2}%
\ell(P(L_{2t+\alpha}))+1$.

If $\ell(P(L_{2t+\alpha}))$ is odd, i.e., if $\alpha$ is odd, then
$\ell(P(H))$ is even; hence $\sigma(h)\geq\frac{1}{2}\ell(P(H))+2$. As in
previous cases, it follows that $L_{2t+\alpha+1}$ and $R_{2t+\alpha+1}$ exist.
Moreover, $\ell(P(L_{2t+\alpha+1}))=\ell(P(L_{2t}))+\ell(P(H))+\ell
(P_{2t+\alpha+1})+2$ and%
\begin{align}
\sigma(g_{2t+\alpha+1})  &  =\sigma(g_{2t})+\sigma(h)+\sigma(f_{2t+\alpha
+1})\nonumber\\
&  \geq\frac{1}{2}[\ell(P(L_{2t}))+\ell(P(H))+\ell(P_{2t+\alpha+1}%
)]+2=\frac{1}{2}\ell(P(L_{2t+\alpha+1}))+1. \label{eq_alphaH}%
\end{align}
Thus, whether or not $P(L_{2t+\alpha+1})$ is a diametrical path of
$L_{2t+\alpha+1}$, we have that\newline$\sigma(g_{2t+\alpha+1})\geq\frac{1}%
{2}\ell(P(L_{2t+\alpha+1}))+1.$

\item Suppose $\alpha=1$. Then $H=T_{2t+1}$. Similar to (\ref{eq_alphaH}), but
only knowing that $\sigma(h)>\frac{1}{2}\ell(P(H))$, we obtain that
\[
\sigma(g_{2t+2})\geq\frac{1}{2}\ell(P(L_{2t+2})).
\]
Suppose $P(L_{2t+2})$ is not a diametrical path of $L_{2t+2}$, and let
$Q(L_{2t+2})$ be such a path instead. Define $z_{2t+2}$ and $q_{2t+2}$ as in
the proof of Claim \ref{Cl_3}. Since $P(L_{2t})$ is a diametrical path of
$L_{2t}$, $z_{2t+2}\in V(T_{2t+1})\cup V(T_{2t+2})$.

\begin{itemize}
\item Suppose $z_{2t+2}\in V(T_{2t+1})$. Let $F$ be the subtree of $T$ induced
by $V(T_{2t+1})\cup V(T_{2t+2})$, with $P(F)$ the subpath of $P$ in $F$, and
$f^{\ast}$ the restriction of $f$ to $F$. As before, $f^{\ast}$ is a maximal
bn-independent broadcast on $F$, hence $\sigma(f^{\ast})\geq i_{\operatorname{bn}%
}(F) $. Note that $\ell(P(F))=\ell
(P_{2t+1})+\ell(P_{2t+2})+1$, which is odd. Then all edges on the
$y_{2t}-z_{2t+2}$ subpath of $P(F)$ are covered, while the edge $x_{2t+1}%
y_{2t+1}$ on the $z_{2t+2}-x_{2t+2}$ subpath of $P(F)$ is uncovered. By
Corollary \ref{Cor_edge_cover2}, $\sigma(f^{\ast})>\frac{1}{2}\ell(P(F))+1$,
that is, $\sigma(f^{\ast})\geq\frac{1}{2}[\ell(P(F))+1]+1$. Now we have
$\ell(P(L_{2t+2}))=\ell(P(L_{2t}))+\ell(P(F))+1$ and%
\[
\sigma(g_{2t+2})=\sigma(g_{2t})+\sigma(f^{\ast})\geq\frac{1}{2}[\ell
(P(L_{2t}))+\ell(P(F))+1]+1=\frac{1}{2}\ell(P(L_{2t+2}))+1.
\]

\item Suppose $z_{2t+2}\in V(T_{2t+2})$. By Corollary \ref{Cor_edge_cover2}%
$(i)$, $\sigma(f_{2t+2})>\frac{1}{2}\ell(P_{2t+2})$. Since $H=T_{2t+1}$, we
also have $\sigma(f_{2t+1})>\frac{1}{2}\ell(P_{2t+1})$. Since $\ell(P_{i})$ is
even for each $i$, we have $\sigma(f_{2t+1})\geq\frac{1}{2}\ell(P_{2t+1})+1$
and $\sigma(f_{2t+2})\geq\frac{1}{2}\ell(P_{2t+2})+1$. Now $\ell
(P(L_{2t+2}))=\ell(P(L_{2t})+\ell(P_{2t+1})+\ell(P_{2t+2})+2$ and
\begin{align*}
\sigma(g_{2t+2})  &  =\sigma(g_{2t})+\sigma(f_{2t+1})+\sigma(f_{2t+2})\\
&  \geq\frac{1}{2}[\ell(P(L_{2t}))+\ell(P_{2t+1})+\ell(P_{2t+2})]+2=\frac
{1}{2}\ell(P(L_{2t+2}))+1.
\end{align*}

\end{itemize}
\end{itemize}

In all cases we have an even integer $2t_{1}>2t$ such that $\sigma(g_{2t_{1}%
})\geq\frac{1}{2}\ell(P(L_{2t_{1}}))$, and if $P(L_{2t_{1}})$ is a not a
diametrical path of $L_{2t_{1}}$, then $\sigma(g_{2t_{1}})\geq\frac{1}{2}%
\ell(P(L_{2t_{1}}))+1$. Thus we can repeat the entire process indefinitely,
constructing a strictly increasing infinite sequence $t,t_{1},...$ which
corresponds to the sequence of subtrees $L_{2t}\subsetneqq L_{2t_{1}%
}\subsetneqq...$ of $T$, which, of course, is impossible. This completes the
proof of Case 2(a).~$\blacklozenge$

\bigskip

\noindent\textbf{Case 2(b):\hspace{0.1in}}Suppose $\ell(P_{1})$ is even. By
Lemma \ref{Lem_edge_cover}$(ii)$, $\sigma(f_{1})\geq\frac{1}{2}\ell(P)$.
Equivalently, considering $L_{1}$, $\sigma(g_{1})\geq\frac{1}{2}\ell
(P(L_{1}))$. If $P(L_{1})$ is not a diametrical path of $L_{1}$, then
$\sigma(g_{1})\geq\frac{1}{2}\ell(P(L_{1}))+1$. The proof now proceeds exactly
as in Case 2(a), except that we consider the subtrees $L_{i}$ and $R_{i}$ for
odd indices $i$ instead of even ones, because $\ell(P(L_{i}))$ is even . The
resulting contradiction concludes the proof of the theorem.~$\blacklozenge
$~$\blacksquare$

\subsection{Comments on Theorems \ref{Thm_ibn_bound} and
\ref{Thm_Uniquely_radial}}

Paths $P_{3k}$, where $k\geq4$, have split-sets of cardinality $k-1$, uniquely
radial subtrees $P_{3}$, and $\gamma_{b}(P_{3k})=k$, but by Theorem
\ref{Thm_paths}, $i_{\operatorname{bn}}(P_{3k})=\left\lceil 6k/5\right\rceil
<\left\lceil 4k/3\right\rceil $. Therefore Theorem \ref{Thm_Uniquely_radial}
cannot be extended to trees with larger split-sets.

The condition that the radial subtrees of $T$ be uniquely radial is used only
in the proof of Claim \ref{Cl_1}. It is not hard to see that the proof remains
valid when we relax this condition to only require that when $m=1$, at least
one of the radial subtrees of $T$ is uniquely radial, while when $m=2$, both
radial subtrees that contain peripheral vertices of $T$ are uniquely radial,
but the subtree that lies between the two split-edges need not be uniquely
radial. For example, for each of the split-sets of the tree in Figure
\ref{Fig_split}, only one of the radial subtrees is uniquely radial; as shown
in Figure~\ref{Fig_bn_ind}, $i_{\operatorname{bn}}(T)=\gamma_{b}%
(T)+\left\lceil 2/3\right\rceil =5$.

\section{Open Problems}

\label{Sec_open}We close with a list of open problems.

\begin{problem}
Find more trees (other than paths) whose radial subtrees are uniquely radial
but for which $i_{\operatorname{bn}}(T)<\gamma_{b}(T)+\left\lceil \frac
{|M|+1}{3}\right\rceil $ when $|M|\geq3$.
\end{problem}

\begin{problem}
Let $T$ be a tree whose radial subtrees are uniquely radial with radii at
least $2$ for any maximum split-set $M$. Is it true that $i_{\operatorname{bn}%
}(T)=\gamma_{b}(T)+\left\lceil \frac{|M|+1}{3}\right\rceil $?
\end{problem}

Herke \cite{Herke} showed that if $G$ is a connected graph, then $\gamma
_{b}(G)=\min\{\gamma_{b}(T):T$ is a spanning tree of $G\}$.

\begin{problem}
Is it true that if $G$ is a connected graph, then $i_{\operatorname{bn}%
}(G)=\min\{i_{\operatorname{bn}}(T):T$ is a spanning tree of $G\}$?
\end{problem}

The proof of the bound $i_{\operatorname{bn}}(G)\leq\left\lceil 4\gamma
_{b}(G)/3\right\rceil $ in \cite{N thesis} shows that equality holds for a
graph $G$ only if $G$ has a $\gamma_{b}$-broadcast $f$ such that $f(v)=1$ for
each $v\in V_{f}^{+}$, $|V_{f}^{+}|\geq2$, and $N_{f}(v)\cap N_{f}%
(u)=\varnothing$ for any distinct $u,v\in V_{f}^{+}$. In this case, $f$ is the
characteristic function of an \emph{efficient dominating set} of $G$. The
paths $P_{6}$ and $P_{9}$ are examples of graphs with efficient dominating
sets for which equality holds in the bound.

\begin{problem}
Determine classes of graphs such that $i_{\operatorname{bn}}(G)=\left\lceil
\frac{4\gamma_{b}(G)}{3}\right\rceil $. Improve the bound for classes of
graphs other than trees.
\end{problem}

As mentioned in the introduction, $i(G)$ and $i_{\operatorname{bn}}(G)$ are
not comparable. Since $i(K_{n,n})=n$ and $i_{\operatorname{bn}}(K_{n,n})=2$
when $n\geq2$, the ratio $i(G)/i_{\operatorname{bn}}(G)$ is, in general,
unbounded. Since $\gamma_{b}(G)\leq\gamma(G)\leq i(G)$ for all graphs $G$, we
also have that $i_{\operatorname{bn}}(G)\leq\left\lceil 4i(G)/3\right\rceil $
for all graphs $G$. Paths show that the difference $i_{\operatorname{bn}%
}(G)-i(G)$ can be arbitrarily large, but $i_{\operatorname{bn}}(P_{n}%
)/i(P_{n})\approx6/5$. In general, the actual value of $i_{\operatorname{bn}%
}(G)$ is likely to be smaller than $4i(G)/3$.

\begin{problem}
Improve the bound $i_{\operatorname{bn}}(G)\leq\left\lceil \frac{4i(G)}%
{3}\right\rceil $. Bound the ratio $i(G)/i_{\operatorname{bn}}(G)$ for
specific graph classes.
\end{problem}

As mentioned in the introduction, Erwin denoted the minimum weight of a
maximal hearing independent broadcast $f$ on $G$ by $i_{b}(G)$. We prefer the
notation $i_{h}(G)$ ($h$ for \textquotedblleft hearing\textquotedblright) to
avoid confusion with other notation and to be consistent with \cite{MN, N
thesis}. A tree $T$ with $i_{\operatorname{bn}}(T)=8$ and $i_{h}%
(T)=\operatorname{rad}(T)=9$ is shown in Figure \ref{Fig_bn_h}; hence these
parameters are distinct. Similar to the result in Theorem \ref{Thm_paths},
Bouchouika, Bouchemakh and Sopena \cite{BBS} showed that $i_{h}(P_{n}%
)=i_{h}(C_{n})=\left\lceil 2n/5\right\rceil $ for all $n\neq3$.

\begin{figure}[ht]
    \centering
    \begin{tikzpicture}[main/.style = {draw, circle}, minimum size =0.25 cm, inner sep=0pt, node distance =0.75cm]
    
    \node[main,fill=blue] (1) {};
    \node[main,fill=blue] (2) [right of =1]{};
    \node[main,fill=blue,very thick] (3) [right of =2,label={[xshift=0.0cm, yshift=-0.75 cm]\textcolor{blue}{2}}]{};
    \node[main,fill=blue] (4) [right of =3]{};
    \node[draw,rectangle,fill=plum] (5) [right of =4]{};
    \node[main,fill=red] (6) [right of =5] {};
    \node[main,fill=red,very thick,label={[xshift=0.0cm, yshift=-0.75 cm]\textcolor{red}{2}}] (7) [right of =6]{};
    \node[main,fill=red] (8) [right of =7]{};
    \node[main,fill=red] (9) [right of =8] {};
    \node[main,fill=blue] (10) [right of =9] {};
    \node[main,fill=blue] (11) [right of =10]{};
    \node[main,fill=blue,very thick] (12) [right of =11,,label={[xshift=0.0cm, yshift=-0.75 cm]\textcolor{blue}{2}}]{};
    \node[main,fill=blue] (13) [right of =12]{};
    \node[draw,rectangle,fill=plum] (14) [right of =13] {};
    \node[main,fill=red] (15) [right of =14]{};
    \node[main,fill=red,very thick] (16) [right of =15,label={[xshift=0.0cm, yshift=-0.75 cm]\textcolor{red}{2}}] {};
    \node[main,fill=red] (17) [right of =16] {};
   \node[main,fill=red] (a) [right  of =17] {};
    
    \node[main,fill=blue] (18) [above of =3] {};
    \node[main,fill=blue] (19) [above of =18] {};
    \node[main,fill=red] (20) [above of =7] {};
    \node[main,fill=red] (21) [above of =20] {};
    \node[main,fill=blue] (22) [above of =12]{};
    \node[main,fill=blue] (23) [above of =22] {};
    \node[main,fill=red] (24) [above of =16] {};
    \node[main,fill=red] (25) [above of =24] {};
    
    \draw(1)--(2);
    \draw(2)--(3);
    \draw(3)--(4);
    \draw(4)--(5);
    \draw(5)--(6);
    \draw(6)--(7);
    \draw(7)--(8);
    \draw(8)--(9);
    \draw[ultra thick](9)--(10);
    \draw(10)--(11);
    \draw(11)--(12);
    \draw(12)--(13);
    \draw(13)--(14);
    \draw(14)--(15);
    \draw(15)--(16);
    \draw(16)--(17);
   {[on background layer] \draw[](16)--(a);}

    \draw(3)--(18);
    \draw(18)--(19);
    \draw(7)--(20);
    \draw(20)--(21);
    \draw(12)--(22);
    \draw(23)--(22);
    \draw(16)--(24);
    \draw(24)--(25);

    \end{tikzpicture}
    \caption{A tree $T$ with $i_{bn}(T)=8$ and $i_h(T)=\mathrm{rad}(T)=9$}
    \label{Fig_bn_h}
\end{figure}

\begin{problem}
$(i)\hspace{0.1in}$For which graphs $G$ is $i_{h}(G)=i_{\operatorname{bn}}(G)$?

\begin{enumerate}
\item[$(ii)$] Can the difference $i_{h}(G)-i_{\operatorname{bn}}(G)$ be arbitrary?\ 

\item[$(iii)$] Is the ratio $i_{h}(G)/i_{\operatorname{bn}}(G)$ bounded?

\item[$(iv)$] Do there exist graphs $G$ such that $i_{\operatorname{bn}%
}(G)>i_{h}(G)$?
\end{enumerate}
\end{problem}

Heggernes and Lokshtanov \cite{HL} showed that minimum broadcast domination is
solvable in polynomial time for any graph. Their algorithm runs in $O(n^{6})$
time for a graph of order $n$.

\begin{problem}
Study the complexity of determining $i_{\operatorname{bn}}(G)$ for $G$
belonging to various graph classes.
\end{problem}

Brewster, Mynhardt and Teshima \cite{BMT} considered broadcast domination as
an integer programming (IP) problem. Its fractional relaxation linear program
(LP) has a dual linear program whose IP formulation provides a lower bound for
the broadcast number.

\begin{problem}
Consider the bn-independent broadcast problem as an IP problem and obtain its
dual IP formulation via relaxation to an LP problem and its dual.
\end{problem}
\smallskip

\noindent\textbf{Acknowledgement\hspace{0.1in}}We acknowledge the support of
the Natural Sciences and Engineering Research Council of Canada (NSERC), PIN 253271.

\noindent Cette recherche a \'{e}t\'{e} financ\'{e}e par le Conseil de
recherches en sciences naturelles et en g\'{e}nie du Canada (CRSNG), PIN
253271.
\begin{center}
\includegraphics[width=2.5cm]{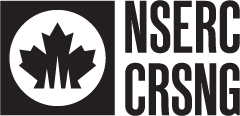}%
\end{center}

\label{refs}

\end{document}